\documentclass[12pt,a4paper,final]{iopart}

\usepackage{color}
\usepackage{iopams}  
\usepackage{graphicx}
\usepackage[breaklinks=true,colorlinks=true,linkcolor=blue,urlcolor=blue,citecolor=blue]{hyperref}

\expandafter\let\csname equation*\endcsname\relax
\expandafter\let\csname endequation*\endcsname\relax

\usepackage{amsmath}
\usepackage{amsthm}
\usepackage[]{algorithm2e}
\usepackage{multirow}
\usepackage[super]{nth}
\usepackage{subcaption}
\usepackage{tabularx}
\usepackage[capitalise]{cleveref}

\usepackage{amssymb}
\usepackage{dsfont}
\usepackage{mathrsfs}
\usepackage{float,indentfirst}
\newtheorem{theorem}{Theorem}
\theoremstyle{plain}

\newtheorem{remark}{Remark}
\newtheorem{definition}{Definition}
\newtheorem{lemma}{Lemma}

\newcommand{\mean}{m}
\newcommand{\pmean}{\widehat{m}}
\newcommand{\pp}{\widehat{\theta}}
\newcommand{\py}{\widehat{y}}

\newcommand{\Cov}{C}
\newcommand{\pCov}{\widehat{C}}

\newcommand{\pD}{\hat{D}}

\newcommand{\CovZ}{Z}
\newcommand{\pCovZ}{\widehat{Z}}

\newcommand{\pGCovZ}{\widehat{\mathcal{Y}}}

\newcommand{\N}{\mathcal{N}}
\newcommand{\G}{\mathcal{G}}

\newcommand{\E}{\mathbb{E}}
\newcommand{\I}{\mathbb{I}}
\newcommand{\R}{\mathbb{R}}
\newcommand{\Z}{\mathbb{Z}}

\newcommand{\bG}{\mathcal{F}}
\newcommand{\blG}{F}

\newcommand{\by}{x}
\newcommand{\bpy}{\widehat{\by}}
\newcommand{\Dt}{\Delta \tau}

\newcommand{\bigO}{\mathcal{O}}

\definecolor{darkred}{rgb}{.7,0,0}
\definecolor{darkblue}{rgb}{0,0,.7}
\definecolor{darkgreen}{rgb}{0,.7,0}
\definecolor{darkbrown}{rgb}{0.8,0.4,0.4}

\newcommand{\dzh}[1]{{\color{black}{#1}}}

\begin{document}

\title{Efficient Derivative-free Bayesian Inference for Large-Scale Inverse Problems}
  
\author{Daniel~Zhengyu~Huang}
\address{California Institute of Technology, Pasadena, CA}
\ead{dzhuang@caltech.edu}
  
\author{Jiaoyang~Huang}
\address{New York University, New York, NY}
\ead{jh4427@nyu.edu}
  
\author{Sebastian Reich}
\address{Universit\"{a}t Potsdam, Potsdam, Germany}
\ead{sebastian.reich@uni-potsdam.de}
  
\author{Andrew M. Stuart}
\address{California Institute of Technology, Pasadena, CA}
\ead{astuart@caltech.edu}

\begin{abstract}

We consider Bayesian inference for large-scale inverse problems, where computational challenges arise from the need for repeated evaluations of an expensive forward model. This renders most Markov chain Monte Carlo approaches infeasible, since they typically require $\bigO(10^4)$ model runs, or more. Moreover, the forward model is often given as a black box or is impractical to differentiate. Therefore derivative-free algorithms are highly desirable.
We propose a framework, which is built on Kalman methodology, to efficiently perform Bayesian inference in such inverse problems. The basic method is based on an approximation of the filtering distribution of a novel mean-field dynamical system, into which the inverse problem is embedded as an observation operator. Theoretical properties are established for linear inverse problems, demonstrating that the desired Bayesian posterior is given by the steady state of the law of the filtering distribution of the mean-field dynamical system, and proving exponential convergence to it. This suggests that, for nonlinear problems which are close to Gaussian, sequentially computing this law provides the basis for efficient iterative methods to approximate the Bayesian posterior. 
Ensemble methods are applied to obtain interacting particle system approximations of the filtering distribution of the mean-field model; and practical strategies to further reduce the computational and memory cost of the methodology are presented, including low-rank approximation and a bi-fidelity approach. The effectiveness of the framework is demonstrated in several numerical experiments, including proof-of-concept linear/nonlinear examples and two large-scale applications: learning of permeability parameters in subsurface flow; and learning subgrid-scale parameters in a global climate model. Moreover, the stochastic ensemble Kalman filter and various ensemble square-root Kalman filters are all employed and are compared numerically. The results demonstrate that the proposed method, based on exponential convergence to the filtering distribution of a mean-field dynamical system, is competitive with pre-existing Kalman-based methods for inverse problems. 
\end{abstract}

%
\vspace{2pc}
\submitto{\IP}
%
%
%

\section{Introduction}

\subsection{Orientation}

The focus of this work is on efficient derivative-free Bayesian inference approaches for large scale inverse problems, in which
the goal is to estimate probability densities for uncertain parameters, given noisy observations derived from the output of
a model that depends on the parameters.
Such approaches are highly desirable for numerous models arising in
science and engineering applications, often defined through
partial differential equations (PDEs).
These include,  to name a few, global climate model calibration~\cite{sen2013global,schneider2017earth}, material constitutive relation calibration~\cite{huang2020learning,xu2021learning,avery2021computationally}, seismic inversion in geophysics~\cite{russell1988introduction,bunks1995multiscale,bui2013computational,martin2012stochastic,wang2022joint}, and biomechanics inverse problems~\cite{bertoglio2012sequential,moireau2011reduced}.
Such problems may feature multiple scales, may include chaotic dynamics, or may involve turbulent phenomena; as a result
the forward models are typically very expensive to evaluate. Moreover, the forward solvers are often given as a black box (e.g., off-the-shelf solvers~\cite{jasak2007openfoam} or multiphysics systems requiring coupling of different solvers~\cite{huang2020modeling,cao2022bayesian}), and may not
be differentiable due to the numerical methods used (e.g., embedded boundary method~\cite{peskin1977numerical,huang2018family} and adaptive mesh refinement~\cite{berger1989local,borker2019mesh}) 
or because of the inherently discontinuous physics (e.g. in fracture~\cite{moes1999finite} or cloud modeling~\cite{tan2018extended,lopez2022training}).

Traditional methods for derivative-free Bayesian inference to estimate the posterior distribution include specific instances
of the Markov chain Monte Carlo methodology~\cite{geyer1992practical,gelman1997weak,goodman2010ensemble,cotter2013mcmc}
(MCMC), such as random walk Metropolis or the preconditioned Crank-Nicolson (pCN) algorithm~\cite{cotter2013mcmc}, and Sequential Monte Carlo methods~\cite{del2006sequential,beskos2015sequential} (SMC),
which are in any case often interwoven with MCMC.
These methods typically require  $\mathcal{O}(10^4)$ iterations, or more, to reach statistical convergence for the complex forward models which
motivate our work.
Given that each forward run can be expensive, conducting $\mathcal{O}(10^4)$ runs is often computationally unfeasible. 
We present an approach based on the Kalman filter methodology, which aims to estimate the first two moments of the posterior distribution. We demonstrate that, in numerical tests across a range of
examples, the proposed methodologies converge within $\bigO(10)$ iterations, using $\bigO(10)$ embarrassingly-parallel model evaluations per step,
resulting in orders of magnitude reduction in cost over derivative-free MCMC and SMC methods.
We also demonstrate favorable performance in comparison with existing Kalman-based Bayesian inversion techniques.

In Subsection \ref{ssec:1.2}, we outline the Bayesian approach
to inverse problems,
describing various approaches to sampling, formulated as dynamical systems
on probability measures, and introducing our novel mean field approach.
In Subsection \ref{ssec:1.3},
we discuss pathwise stochastic dynamical systems
which realize such dynamics at the level of measures, and discuss filtering
algorithms which may be applied to them for the purposes of approximate
inversion.
Subsection \ref{ssec:1.4} highlights the novel contributions
in this paper, building on the context established in the two
preceding subsections. Subsection \ref{ssec:1.5} summarizes notational
conventions that we adopt throughout.

\subsection{Bayesian Formulation of the Inverse Problem}
\label{ssec:1.2}

Inverse problems can be formulated as recovering unknown parameters $\theta \in \R^{N_{\theta}}$ from noisy observation $y \in \R^{N_y}$ related through
\begin{equation}
\label{eq:KI}
    y = \G(\theta) + \eta.
\end{equation}
Here $\G$ denotes a forward model mapping parameters to output observables, and $\eta$ denotes observational noise; for
simplicity we will assume known Gaussian statistics: $\eta \sim \N(0,\Sigma_{\eta})$.
In the Bayesian perspective, $\theta$ and $y$ are treated as random variables. Given the prior $\rho_{\rm prior}(\theta)$ on $\theta$, the inverse problem can be formulated as finding the posterior $\rho_{\rm post}(\theta)$ on $\theta$ given $y$~\cite{kaipio2006statistical, dashti2013bayesian,dashti2013map}:
\begin{equation}
\label{eq:KI3}
\rho_{\rm post}(\theta) = \frac{1}{Z(y)} e^{-\Phi(\theta , y) } \rho_{\rm prior}(\theta), \qquad  \Phi(\theta,y) = \frac{1}{2}\lVert\Sigma_{\eta}^{-\frac{1}{2}}(y - \G(\theta)) \rVert^2
\end{equation}
and $Z(y)$ is the normalization constant
\begin{equation}
    Z(y) =  \int e^{-\Phi(\theta , y) } \rho_{\rm prior}(\theta)d\theta.
\end{equation}
We focus on the case, where the prior $\rho_{\rm prior}$ is (or is approximated as) Gaussian with  mean and covariance
$r_0$ and $\Sigma_{0}$, respectively. Then the posterior $\rho_{\rm post}(\theta)$  can be written as 
\begin{equation}
\label{eq:KI2}
    \rho_{\rm post}(\theta) = \frac{1}{Z(y)} e^{-\Phi_R(\theta , y) }, \qquad  \Phi_R(\theta,y) = \Phi(\theta,y)+\frac{1}{2}\lVert\Sigma_{0}^{-\frac{1}{2}}(\theta - r_0) \rVert^2.
\end{equation}

\subsubsection{Computational Approaches}
Bayesian inference requires approximation of, or 
samples from, the posterior distribution given by~\cref{eq:KI3}. There are three major avenues to approximate the posterior distribution: 
\begin{itemize}
    \item Those based on variational inference~\cite{anderson1987mean,parisi1988statistical}, where a parameterized approximate density is constructed and optimized to minimize the distance to the posterior density. They include Gaussian variational inference~\cite{opper2009variational,quiroz2018gaussian,galy2021particle} and normalizing  flows~\cite{rezende2015variational}.
    
    \item Those based on sampling and more importantly the invariance of measures and ergodicity. They include MCMC~\cite{geyer1992practical,gelman1997weak}, Langevin dynamics~\cite{rossky1978brownian,roberts1996exponential}, and more recently interacting particle approaches~\cite{vrugt2009accelerating,goodman2010ensemble,foreman2013emcee,LMW18}.
    
    At an abstract mathematical level, invariance and ergodicity-based approaches to sampling from the posterior $\rho_{\rm post}$ rely on the transition kernel $\psi_{\rm I}(\theta', \theta)$ such that
\begin{equation} \label{eq:invariance}
    \rho_{\rm post}(\theta) = \int \psi_{\rm I}(\theta', \theta) \rho_{\rm post}(\theta') 
    d\theta',
\end{equation}
that is, the posterior distribution $\rho_{\rm post}(\theta)$ is  
invariant with respect to the transition kernel $\psi_{\rm I}(\theta',\theta)$.
Furthermore, starting from any initial distribution the associated Markov chain should approach the invariant measure $\rho_{\rm post}(\theta)$ asymptotically. 

    \item Those based on coupling ideas (mostly in the form of coupling the prior with the posterior). While several sequential data assimilation methods, such as importance sampling-resampling in SMC~\cite{smith2013sequential} and the ensemble Kalman filtering~\cite{evensen2009ensemble,reich2015probabilistic,law2015data}, can be viewed under the coupling umbrella, the systematic exploitation/exposition of the coupling perspective in the context of Bayesian inference is more recent, including the ideas of transport maps~\cite{reich2011dynamical,el2012bayesian,reich2013nonparametric,marzouk2016sampling,ruchi2019transform}.
    
    At an abstract mathematical level, the coupling approach is based on a transition kernel $\psi_{\rm C}(\theta',\theta)$ such that
\begin{equation} \label{eq:coupling}
    \rho_{\rm post}(\theta) = \int \psi_{\rm C}(\theta',\theta)
    \rho_{\rm prior}(\theta')d\theta'.
\end{equation}
The transition kernel forms a coupling between the prior and the posterior distribution and is applied only once. The induced transition from $\theta' \sim \rho_{\rm prior}$ to $\theta \sim \rho_{\rm post}$ is of the type of a McKean--Vlasov mean-field process and can be either deterministic or stochastic \cite{reich2015probabilistic}.
In practice the methodology is implemented via an approximate coupling, using linear transport maps:
\begin{equation} \label{eq:linear_coupling}
    \theta = A\theta' + b ,
\end{equation}
where the matrix $A$ and the vector $b$ depend on the prior distribution $\rho_{\rm prior}$, the data likelihood $\Phi$, and the data $y$, and are chosen such that the induced random variable $\theta$ approximately samples from the posterior distribution $\rho_{\rm post}$. Many variants of the popular ensemble Kalman filter can be derived within this framework. 

\end{itemize}

\subsubsection{A Novel Algorithmic Approach}
\label{sec:novel_approach}

The main contribution of this paper is to incorporate all three approaches from above by designing a particular (artificial) mean-field dynamical system and applying filtering methods, which employ
a Gaussian ansatz, to approximate the filtering distribution resulting from partial observation of the system; the equilibrium of the
filtering distribution is designed to be close to the desired
posterior distribution. At an abstract level, we introduce a data-independent transition kernel, denoted by $\psi_{\rm P}(\theta'', \theta')$, and another data-dependent transition kernel, denoted by $\psi_{\rm A}(\theta'', \theta')$, such that the posterior distribution $\rho_{\rm post}$ remains invariant under the both transition kernels combined, that is,
\begin{equation} \label{eq:KI_invariance}
    \rho_{\rm post}(\theta) =
    \int  \psi_{\rm A}(\theta', \theta) \left(\int
    \psi_{\rm P}(\theta'', \theta') \rho_{\rm post}(\theta'')d\theta''\right) d \theta'.
\end{equation}
The first transition kernel, 
$\psi_{\rm P}(\theta'', \theta')$, corresponds to the prediction step in filtering methods and is chosen such that
\begin{equation}
\label{eq:KI_pred_0}
    \hat \rho_{n+1}(\theta) = \int \psi_{\rm P}(\theta',\theta)\rho_n (\theta')d \theta' \quad \mbox{and} \quad \hat \rho_{n+1}(\theta) \propto \rho_n (\theta)^{\dzh{ 1-\Dt}} ,
\end{equation}
where \dzh{$0 < \Dt < 1$ is the time-step size, a free parameter,} and $\rho_n(\theta)$ denotes the current density. In other words, this transition kernel corresponds to a simple rescaling of a given density. The second transition kernel, $\psi_{\rm A}(\theta'',\theta')$, corresponds to the analysis step in filtering methods and has to satisfy
\begin{equation} 
\label{eq:KI_analysis_0}
\rho_{n+1}(\theta) = \int \psi_{\rm A}(\theta',\theta)\hat \rho_{n+1} (\theta')d\theta'
\quad \mbox{and} \quad \rho_{n+1}(\theta) \propto \rho_{\rm post}(\theta)^{\dzh{\Dt}}\hat \rho_{n+1}(\theta).
\end{equation}
This transition kernel depends on the data and the posterior distribution and performs a suitably modified Bayesian inference step. Combining the
two preceding displays yields
\begin{align}
\rho_{n+1}(\theta) \propto \rho_{\rm post}(\theta)^{\dzh{\Dt}} \rho_{n}(\theta)^{\dzh{1 - \Dt}}.
\end{align}
It is immediate that the overall transition $\rho_n \mapsto \rho_{n+1}$ is
indeed invariant with respect to $\rho_{\rm post}$; furthermore,
by taking logarithms in the mapping
from $\rho_{n}$ to $\rho_{n+1}$ it is possible to deduce exponential convergence to this steady state, for any $\dzh{0 < \Dt < 1}$. In our concrete algorithm a mean field dynamical system is introduced for which \cref{eq:KI_analysis_0} is satisfied exactly, whilst  \cref{eq:KI_pred_0}
is satisfied only in the linear, Gaussian setting; the resulting
filtering distribution is approximated using Kalman methodology
applied to filter the resulting partially observed mean-field dynamical system. We emphasize that the involved transition kernels are all of McKean--Vlasov type, that is, they depend on the distribution of the parameters $\theta$.

There are several related approaches. We mention in this context in particular the recently proposed consensus-based methods. These sampling methods were analyzed in the context of optimization in \cite{carrillo2018analytical}. Similar ideas were then developed for consensus based sampling (CBS) \cite{carrillo2021consensus_sampling} based on the same principles employed here: to find a mean-field model
which, in the linear Gaussian setting converges asymptotically to the posterior distribution, and then to develop implementable algorithms by employing finite
particle approximations of the mean-field. Another related approach has been proposed in \cite{pidstrigach2021affine} where data assimilation algorithms are combined with stochastic dynamics in order to approximately sample from the posterior distribution $\rho_{\rm post}$.

\subsection{Filtering Methods for Inversion}
\label{ssec:1.3}

Since filtering methods are at the heart of our proposed methodology, we provide here a brief summary of a few key concepts. Filtering methods may be deployed to approximate the posterior distribution given by~\cref{eq:KI3}.
The inverse problem is first paired with a dynamical system for the parameter~\cite{emerick2013investigation,chen2012ensemble,iglesias2013ensemble,wan2000unscented}, leading to a hidden Markov
model, to which filtering methods may be applied. In its most basic form, the hidden Markov model takes the form
\begin{subequations}
\label{eq:KL_dynamics}
  \begin{align}
  &\textrm{evolution:}    &&\theta_{n+1} = \theta_{n}, \\
  &\textrm{observation:}  &&y_{n+1} = \G(\theta_{n+1}) + \eta_{n+1}; 
\end{align}
\end{subequations}
here $\theta_{n}$ is the unknown state vector, \dzh{$y_{n+1}$ is the output of the observation model}, and $\eta_{n+1} \sim \N(0, \Sigma_{\eta})$ is the observation error at the $n^{th}$ iteration. \dzh{Any filtering method can be applied to estimate $\theta_n$ given observation data $\{y_{\ell}^{\dagger}\}_{\ell=1}$.} The Kalman filter~\cite{kalman1960new} can be applied to this setting provided the forward operator $\G$ is linear and the initial state $\theta_0$ and the observation errors are Gaussian. The Kalman filter has been extended to nonlinear and non-Gaussian settings in manifold ways,
including but not limited to, the extended Kalman filter (EKF, or sometimes ExKF)~\cite{sorenson1985kalman,jazwinski2007stochastic}, the ensemble Kalman filters (EnKF) ~\cite{evensen1994sequential,anderson2001ensemble,bishop2001adaptive}, and the unscented Kalman (UKF) filter~\cite{julier1995new,wan2000unscented}. We refer to the extended, ensemble and unscented Kalman filters as {\em approximate Kalman filters} to highlight the fact that, outside the linear setting where the Kalman filter \cite{kalman1960new} is exact, they are all uncontrolled approximations designed on the principle of matching first and second moments. 

More precisely, the EnKF uses Monte Carlo sampling to estimate desired means and covariances empirically. Its update step is of the form (\ref{eq:linear_coupling}) and can be either deterministic or stochastic. The ensemble adjustment/transform filters are particle approximations of square root filters, a deterministic approach to matching first and second moment information \cite{tippett2003ensemble}.
The unscented Kalman filter uses quadrature, and is also a deterministic method; it may also be viewed as approximating a square root filter. The stochastic EnKF on the other hand compares the data $y$ to model generated data and its update step is intrinsically  stochastic, that is, the vector $b$ in (\ref{eq:linear_coupling}) itself is random. 

All of the filtering methods to estimate $\theta_n$ given \dzh{$\{y_\ell^{\dagger}\}_{\ell=1}^n$ that we have described so far
may be employed in the setting where $y_\ell^{\dagger} \equiv y$}; repeated exposure of the parameter to the data helps the system to learn the parameter
from the data. In order to maintain statistical consistency, an $N$-fold insertion of the same data $y$ requires an appropriate modification of the data likelihood function and the resulting Bayesian inference step becomes 
\begin{equation}
    \label{eq:Ifilter}
\rho_n(\theta) \rightarrow \rho_{n+1}(\theta)\propto \rho_n(\theta) e^{-\frac{1}{N}\Phi(\theta,y)}.
\end{equation}
Initializing with $\rho_0(\theta) = \rho_{\rm prior}(\theta)$, after $N$ iterations, $\rho_N(\theta)$ is equal to the posterior density. The
filtering distribution for \eqref{eq:KL_dynamics} recovers this exactly if 
\dzh{$y_\ell^\dagger \equiv y$} and if the variance of $\eta$ is rescaled by $N;$ use of ensemble Kalman methods
in this setting leads to approximate Bayesian inference, which is
intuitively accurate when the posterior is close to Gaussian.
We note that the resulting methodology can be viewed as a homotopy method, such as SMC \cite{del2006sequential} and transport
variants \cite{reich2011dynamical}, which seek to deform the prior into the posterior in one unit time with a finite number of inner steps 
$N$ -- foundational papers introducing ensemble
Kalman methods in this context are \cite{emerick2013investigation,chen2012ensemble,bocquet2013joint}. 
\dzh{Adaptive time-stepping strategies in this context are explored in
\cite{iglesias2021adaptive,matveev2021bayesian,tso2021efficient}.}
Throughout this paper, we will denote the resulting methods as iterative extended Kalman filter (IEKF), iterative ensemble Kalman filter (IEnKF), iterative unscented Kalman filter (IUKF), iterative ensemble adjustment Kalman filter (IEAKF)
and  iterative ensemble transport Kalman filter (IETKF).

We emphasize that multiple insertions of the same data $y$ without the adjustment (\ref{eq:Ifilter}) of the data likelihood function, and/or over
arbitrary numbers of steps, leads to the class of optimization-based Kalman inversion methods: EKI~\cite{iglesias2013ensemble}, Tikhonov-regularized EKI, termed TEKI~\cite{chada2020tikhonov} and unscented Kalman inversion, UKI~\cite{UKI}; see  also \cite{weissmann2022adaptive} for recent adaptive
methodologies which are variants on TEKI. These variants of the Kalman filter lead to efficient derivative-free optimization approaches to approximating the maximum likelihood estimator or maximum a posteriori estimator in the asymptotic limit as $n \to\infty$.
The purpose of our paper is to develop similar ideas, based on
iteration to infinity in $n$, but to tackle the problem of sampling
from the posterior $\rho_{\rm post}(\theta)$ rather than the optimization problem. To achieve these we introduce a novel mean-field
stochastic dynamical system, generalizing \eqref{eq:KL_dynamics} and
apply ensemble Kalman methods to it. This leads to
Bayesian analogues of EKI and the UKI. To avoid proliferation of
nomenclature, we will also refer to these as EKI and UKI relying on the
context to determine whether the optimization or Bayesian approach is
being adopted; in this paper our
focus is entirely on the Bayesian context. 
We will also use ensemble adjustment and transform filters,
denoted as EAKF and ETKF, noting that these two may be applied
in either the optimization (using \eqref{eq:KL_dynamics}) or Bayesian
(using the novel mean-field
stochastic dynamical system introduced here) context, 
but that here we only study the Bayesian problem.
The main conclusions of our work are two-fold, concerning the application of
Kalman methods to solve the Bayesian inverse problem: that with carefully chosen
underlying mean-field dynamical system, such that the prediction and analysis steps approximate~\cref{eq:KI_pred_0}
and replicate ~\cref{eq:KI_analysis_0},
iterating to infinity leads to more efficient and robust methods 
than the homotopy methods
which transport prior to posterior in a finite number of steps; and that deterministic implementations of ensemble Kalman methods, and variants, are superior to stochastic methods. 

The methods we propose are exact in the setting of linear $\G$ and Gaussian prior density $\rho_{\rm prior}$; but,
for nonlinear $\G$, the Kalman-based filters we employ
generally do not converge to the exact posterior distribution, due to the Gaussian ansatz used when deriving the method; negative theoretical results and numerical evidence are reported in~\cite{ernst2015analysis,garbuno2020interacting}. Nonetheless, practical experience demonstrates that the methodology can be effective for problems with distributions close to Gaussian,
a situation which arises in many applications. 

Finally, we note that we also include comparisons with the ensemble Kalman sampler
\cite{garbuno2020interacting,garbuno2020affine,nusken2019note}, which we refer to as the EKS,
an ensemble based Bayesian inversion method derived from discretizing a mean-field stochastic differential equation and which is also based on iteration to infinity, that is, on the invariance principle of the posterior distribution; and we include comparison with the CBS approach \cite{carrillo2021consensus_sampling} mentioned above, another methodology
which also iterates a mean-field dynamical system to infinity to
approximate the posterior.

\subsection{Our Contributions}
\label{ssec:1.4}

The key idea underlying this work is the development of an efficient derivative-free Bayesian inference approach based on applying
Kalman-based filtering methods to a hidden Markov model
arising from a novel mean-field dynamical system.
Stemming from this, our main contributions are as follows:
\footnote{In making these statements, we acknowledge that for linear Gaussian
problems it is possible to solve the Bayesian inverse problem 
exactly in one step,
or multiple steps, using the Kalman filter in transport/coupling mode,
when initialized correctly and with a large enough ensemble. However, the transport/coupling methods are not robust to perturbations from initialization, non-Gaussianity and so forth, whereas the methods we introduce are. Our results substantiate this claim.}

\begin{enumerate}

    \item In the setting of linear Gaussian inverse problems, we prove that the filtering distribution of the mean field model converges exponentially fast to the posterior distribution.
    
    \item We generalize the inversion methods EKI, UKI, EAKI and ETKI
    from the optimization to the Bayesian context by applying the relevant
    variants on Kalman methodologies to the novel mean-field dynamical system
    (Bayesian) rather than to \eqref{eq:KL_dynamics} (optimization).
    
    \item We study and compare application of
    both deterministic and stochastic Kalman methods to the novel mean-field dynamical system, demonstrating that the deterministic methods (UKI, EAKI and ETKI)
    outperform the stochastic method (EKI); this may be attributed to smooth, noise-free approximations resulting from deterministic approaches.
    
    \item We demonstrate that the application of Kalman methods to the
    novel mean-field dynamical system outperforms the application of 
    Kalman filters to transport/coupling models -- the
    IEnKF, IUKF, IEAKF and IETKF approaches; this may be
    attributed to the exponential convergence underlying the
    filter for the novel mean-field dynamical system.
    
    \item We also demonstrate that the application of Kalman methods to the
    novel mean-field dynamical system outperforms the EKS, when Euler-Maruyama discretization is used, because
    the continuous-time formulation requires very small time-steps, and
    CBS which suffers from stochasticity, similarly to the EKI.

    \item We propose several strategies, including low-rank approximation and a bi-fidelity approach, to reduce the computational and memory cost.
      
    \item We demonstrate, on both linear and nonlinear model problems
    (including inference for subsurface geophysical properties in porous
    medium flow), that
    application of deterministic Kalman methods to approximate the
    filtering distribution of the
    novel mean-field dynamical system delivers mean and covariance which are
    close to the truth or to those obtained with the pCN MCMC method.
    The latter uses $\bigO(10^4)$
    model evaluations or more whilst for our method
    only $\bigO(10)$ iterations are required with 
$\bigO(10)$ ensemble members, leading to only $\bigO(10^2)$ 
model evaluations, two orders of magnitude savings.

    \item The method is applied to perform
    Bayesian  parameter inference of
    subgrid-scale parameters arising in an idealized global climate model,
    a problem currently far beyond the reach of state-of-the-art MCMC
    methods such as pCN and variants.
\end{enumerate}

The remainder of the paper is organized as follows. 
In \cref{sec:Alg}, 
the mean field dynamical system, various algorithms which
approximate its filtering distribution, and a complete analysis in the
linear setting, are all presented. \dzh{These correspond to our contributions (i) and (ii). }
In \cref{sec:Variance}, strategies to speed up the algorithm and improve the robustness for real-world problems are presented.
 \dzh{These correspond to our contribution (vi). }
Numerical experiments are provided in \cref{sec:app}; these
serve to empirically confirm the theory and demonstrate the effectiveness of the framework for Bayesian inference. 
\dzh{These correspond to our contributions (iii), (iv), (v), (vii) and (viii).}
We make concluding remarks in
\cref{sec:C}.

The code is accessible online:
\begin{center}
  \url{https://github.com/Zhengyu-Huang/InverseProblems.jl}
\end{center}

\subsection{Notational Conventions}
\label{ssec:1.5}

$A\succ B$ and $A\succeq B$ denote $A-B$ positive-definite or positive-semidefinite, for symmetric matrices $A,B.$ 
$\|\cdot\|, \langle \cdot, \cdot \rangle$ denote Euclidean norm
and inner-product. We use $\mathbb{Z}^+=\{0,1,2,\cdots\}$ to
denote the set of natural numbers; $\N(\cdot, \cdot)$ to denote Gaussian distributions; and $\varrho(\cdot)$ to denote the spectral radius. 
As encountered in Subsection \ref{ssec:1.2} we make  use of the similar symbol $\rho$ for densities; these should be easily distinguished from spectral radius by context and by a different font.

\section{Novel Algorithmic Methodology}

Our novel algorithmic methodology is introduced in this section. 
We first introduce the underlying mean-field dynamical system, which has prediction and analysis steps corresponding to the aforementioned transition kernels, in
Subsection \ref{ssec:MFD}. Then, in Subsection \ref{sec:KF_Intro}, we introduce a class of conceptual Gaussian approximation algorithms 
found by applying Kalman methodology to the proposed mean-field dynamical system. Through linear analysis, we prove in Subsection \ref{ssec:lin} that these algorithms converge exponentially to the posterior. For the
nonlinear setting, a variety of nonlinear Kalman inversion methodologies are discussed in Subsection \ref{sec:nonlinear_Kalman}.

\label{sec:Alg}
\subsection{Mean-Field Dynamical System}
\label{ssec:MFD}

Following the discussion from Section \ref{sec:novel_approach}, we propose an implementation of (\ref{eq:KI_invariance}) to solve inverse problems by pairing the parameter-to-data map with a dynamical system for the parameter, and then employ techniques from filtering to estimate the parameter given the data.

We introduce the prediction step
\begin{equation}
  \theta_{n+1} = \theta_{n} +  \omega_{n+1}. \label{eq:evolve}
\end{equation}
Here $\theta_{n+1}$ is the unknown state vector and $\omega_{n+1} \sim \N(0, \Sigma_{\omega,n+1})$ is the independent, zero-mean Gaussian evolution error, which will be chosen such that (\ref{eq:evolve}) mimics (\ref{eq:KI_pred_0}) for Gaussian densities. The analysis step (\ref{eq:KI_analysis_0}) follows exactly from the  observation model (first introduced in \cite{chada2020tikhonov})
\begin{equation}
    \by_{n+1} = \bG(\theta_{n+1}) + \nu_{n+1}; \label{eq:obser}
\end{equation}
here we have defined the augmented forward map
\begin{align}
    \bG(\theta) = \begin{bmatrix}
        \G(\theta)\\
        \theta
        \end{bmatrix},
\end{align}
with $\nu_{n+1} \sim \N(0, \Sigma_{\nu,n+1})$ the independent, zero-mean Gaussian observation error
\dzh{and $x_{n+1}$ the output of the observation model at time $n+1$. We define artificial observation $\by_{n+1}^{\dagger}$ using the following
particular instance of the data, constructed from the one
observation $y$ and the prior mean $r_0$ and assumed to hold
for all $n \ge 1:$
\begin{align}
\label{eq:data}
    \by_{n+1}^{\dagger} =  \by := \begin{bmatrix}
        y\\
        r_0
        \end{bmatrix}.
    \end{align}
We will apply
filtering methods to condition $\theta_{n}$ on
$Y_{n} := \{\by_1^{\dagger}, \by_2^{\dagger},\cdots , \by_{n}^{\dagger}\}$, the observation set at time $n$. As we will see later, the choice of $\{\by_{\ell}^{\dagger}\}_{l=1}$ leads to the correct posterior.
}

Let $C_n$ denote the covariance of the conditional
    random variable $\theta_n|Y_{n}$. Then the error covariance matrices $\{\Sigma_{\omega,n+1}\}$ and $\{\Sigma_{\nu,n+1}\}$ in the extended dynamical system (\ref{eq:evolve})-(\ref{eq:obser}) are chosen at the $n^{th}$ iteration, as follows:
\begin{equation}
\label{eq:hyperparameters}
    \Sigma_{\nu, n+1}
    = \dzh{\frac{1}{\Dt}}\begin{bmatrix}
        \Sigma_{\eta} & 0\\
        0 & \Sigma_{0} 
        \end{bmatrix}
         \quad \textrm{ and } \quad \Sigma_{\omega,n+1} =  \dzh{\frac{\Dt}{1 - \Dt}}\Cov_{n}.
\end{equation}
Here $\dzh{0 < \Dt < 1}$, and in our numerical studies
we choose $\dzh{\Dt = 1/2}$, although other choices are possible.
Since the artificial evolution error covariance~$\Sigma_{\omega,n+1}$ in \eqref{eq:evolve} 
is updated based on \dzh{$\Cov_{n}$, the conditional covariance of $\theta_n|Y_n$,} it follows
that \eqref{eq:evolve} is a mean-field
dynamical system: it depends on its own law,
specifically on the law of $\theta_n|Y_{n}$.
Details underpinning the choices of the error covariance matrices $\Big\{\Sigma_{\omega,n}\Big\}$ and $\Big\{\Sigma_{\nu,n}\Big\}$ are
given in Subsections~\ref{sec:KF_Intro} and~\ref{ssec:lin}:
the matrices are chosen so that, for linear Gaussian problems,
the prediction and analysis steps follow~\cref{eq:KI_pred_0,eq:KI_analysis_0}, and the converged mean and 
covariance of the resulting filtering distribution for
$\theta_n|Y_{n}$ under the prediction step \eqref{eq:evolve} and the observation model \eqref{eq:obser} match the posterior mean and covariance.

\subsection{Gaussian Approximation}
\label{sec:KF_Intro}
Denote by $\rho_{n}$, the conditional density of $\theta_{n}|Y_{n}$.
We first introduce a class of conceptual Kalman inversion algorithms which approximate $\rho_{n}$ by considering only first and
second order statistics (mean and covariance), and update
$\rho_n$ sequentially using the standard prediction and analysis steps~\cite{reich2015probabilistic,law2015data}:
$\rho_n \mapsto \hat{\rho}_{n+1}$, and then $\hat{\rho}_{n+1} \mapsto \rho_{n+1}$, where $\hat{\rho}_{n+1}$ is the distribution of $\theta_{n+1}|Y_n$. The second analysis step is performed
by invoking a Gaussian hypothesis. In subsequent subsections,
we then apply different methods to
approximate the resulting maps on measures, leading to unscented, stochastic ensemble Kalman and adjustment/transform square root Kalman filters.

In the prediction step, 
assume that $\rho_n \approx \N(\mean_n, \Cov_n)$, then under~\cref{eq:evolve}, $\hat{\rho}_{n+1} = \N(\pmean_{n+1}, \pCov_{n+1})$ is also Gaussian and satisfies
\begin{equation}
\label{eq:KF_pred}
\begin{split}
&\pmean_{n+1} = \E[\theta_{n+1}|Y_n] =  \mean_n \qquad \pCov_{n+1} = \mathrm{Cov}[\theta_{n+1}|Y_n] = \Cov_{n} + \Sigma_{\omega,n+1}.
\end{split}
\end{equation}

In the analysis step, we assume that the joint distribution of     $\{\theta_{n+1}, \by_{n+1}\} | Y_{n}$ can be approximated by a Gaussian distribution
\begin{equation}
\label{eq:KF_joint}
    \N\Big(
    \begin{bmatrix}
    \pmean_{n+1}\\
    \bpy_{n+1}
    \end{bmatrix}, 
    \begin{bmatrix}
  \pCov_{n+1} & \pCov_{n+1}^{\theta \by}\\
    {{\pCov_{n+1}}^{\theta \by}}{}^{T} & \pCov_{n+1}^{\by\by}
    \end{bmatrix}
    \Big),
\end{equation}
where  
\begin{equation}
\label{eq:KF_joint2}
\begin{split}
    &\bpy_{n+1} = \E[\by_{n+1}|Y_n] = \E[\bG(\theta_{n+1})|Y_n], \\
    &\pCov_{n+1}^{\theta \by} = \mathrm{Cov}[\theta_{n+1}, \by_{n+1}|Y_n] = \mathrm{Cov}[\theta_{n+1}, \bG(\theta_{n+1})|Y_n], \\
    &\pCov_{n+1}^{\by\by} = \mathrm{Cov}[\by_{n+1}|Y_n] = \mathrm{Cov}[\bG(\theta_{n+1})|Y_n] + \Sigma_{\nu,n+1}.\\
\end{split}
\end{equation}
These expectations are computed by assuming $\theta_{n+1}|Y_n \sim\hat{\rho}_{n+1}$ and noting that the distribution of $(\theta_{n+1},\by_{n+1})$ is 
then defined by \eqref{eq:evolve}-\eqref{eq:obser}.
This corresponds to projecting
\footnote{We use the term ``projecting'' as finding the Gaussian
$p$ which matches the first and second moments of a given 
measure $\pi$ corresponds
to finding the closest Gaussian $p$ to $\pi$ 
with respect to variation in the
second argument of the (nonsymmetric) Kullback-Leibler divergence
\cite{sanz2018inverse}[Theorem 4.5].}
the joint distribution onto
the Gaussian which matches its mean and covariance. Conditioning the Gaussian in \cref{eq:KF_joint} to find $\theta_{n+1}|\{Y_n, \by_{n+1}^{\dagger}\} = \theta_{n+1}|Y_{n+1}$, gives the following expressions for the mean $\mean_{n+1}$ and covariance $\Cov_{n+1}$ of the approximation to $\rho_{n+1}$ : 
\begin{subequations}
\label{eq:sds}
\begin{align}
        \mean_{n+1} &= \pmean_{n+1} + \pCov_{n+1}^{\theta \by} (\pCov_{n+1}^{\by \by})^{-1} (\by_{n+1}^\dagger - \bpy_{n+1}),                   \label{eq:KF_analysis_mean}\\
         \Cov_{n+1} &= \pCov_{n+1} - \pCov_{n+1}^{\theta \by}(\pCov_{n+1}^{\by \by})^{-1} {\pCov_{n+1}^{\theta \by}}{}^{T}. \label{eq:KF_analysis_cov}
\end{align}
\end{subequations}
\Cref{eq:KF_pred,eq:KF_joint,eq:KF_joint2,eq:KF_analysis_mean,eq:KF_analysis_cov} establish a class of conceptual algorithms for application of Gaussian approximation to solve the inverse problems.
To make implementable algorithms a high level choice needs to be
made: whether to work strictly within the class of
Gaussians, that is to impose $\rho_n \equiv \N(\mean_n, \Cov_n)$, or
whether to allow non-Gaussian $\rho_n$ but to insist that the
second order statistics of the resulting measures agree with
\Cref{eq:KF_pred,eq:KF_joint,eq:KF_joint2,eq:KF_analysis_mean,eq:KF_analysis_cov}. In what follows the UKI takes the first perspective;
all other methods take the second perspective. For the UKI the
method views \Cref{eq:KF_pred,eq:KF_joint,eq:KF_joint2,eq:KF_analysis_mean,eq:KF_analysis_cov} as providing a nonlinear map 
$(m_n,C_n) \mapsto (m_{n+1},C_{n+1})$; this map is then approximated
using quadrature. For the remaining methods a mean-field dynamical
system is used, which is non-Gaussian but matches the aforementioned
Gaussian statistics; this mean-field model is then approximated by a finite
particle system \cite{calvello22}. The dynamical system is of
mean-field type because of the expectations required to calculate
\Cref{eq:KF_joint,eq:KF_joint2} and \eqref{eq:hyperparameters}.
The continuous time limit of the evolution for the mean and covariance
is presented in~\ref{appendix:CTL}; this is obtained by letting $\Dt \to 0.$

\begin{remark}
Consider the case, where $\rho_n = \N(\mean_n, \Cov_n)$ is Gaussian.
With the specific choice of $\Big\{\Sigma_{\omega,n}\Big\}$, we have $\hat{\rho}_{n+1} = \N(\mean_n, \frac{1}{1 - \Dt}\Cov_n)$ from the prediction step~\cref{eq:KF_pred}, and hence the Gaussian density functions $\rho_n$ and  $\hat{\rho}_{n+1}$ fulfill \cref{eq:KI_pred_0}.
With the extended observation model~\eqref{eq:obser} and  the specific choice of  $\Big\{\Sigma_{\nu,n}\Big\}$, the analysis step without Gaussian approximation can be written as 
\begin{align}
\begin{split}
    \rho(\theta_{n+1}|Y_{n+1}) &\propto \rho(\theta_{n+1} | Y_n) \rho(x_{n+1}^\dagger|\theta_{n+1}, Y_n) \\
                               &\propto  \hat{\rho}_{n+1}(\theta_{n+1})  e^{-\Dt\Phi_R(\theta_{n+1}, y)} \\
                               &\propto  \hat{\rho}_{n+1}(\theta_{n+1})  \rho(\theta_{n+1})^{\Dt}
\end{split}
\end{align}
and hence the density functions $\hat{\rho}_{n+1}$ and $\rho_{n+1}$ fulfill \cref{eq:KI_analysis_0}. Note, however, that $\rho_{n+1}$ is not, in general, Gaussian, unless $\G$ is linear, the case studied in the next section.
In the nonlinear case, we employ Kalman-based methodology which only
employs first and second order statistics, and in effect projects
$\rho_{n+1}$ onto a Gaussian.

\end{remark}


\subsection{Linear Analysis}
\label{ssec:lin}
In this subsection, we study the algorithm in the context of linear inverse problems, for which $\G(\theta) = G\theta$ for some
matrix $G.$ Furthermore we assume that $\rho_{\rm prior}$ is Gaussian
$\N(r_0,\Sigma_{0})$ and  recall that the observational noise
is $\N(0,\Sigma_{\eta}).$
Thanks to the linear Gaussian structure the posterior is also
Gaussian with mean and precisions given by 
\begin{equation}
\label{eq:linear_posterior}
   m_{\rm post} = r_0 + \Big(G^T \Sigma_{\eta}^{-1}  G + \Sigma_0^{-1}\Big)^{-1} G^T\Sigma_{\eta}^{-1}(y -  G r_0)\quad \textrm{ and } \quad
   C_{\rm post}^{-1} = G^T\Sigma_{\eta}^{-1}G + \Sigma_0^{-1}
\end{equation}

Furthermore, 
the equations~\eqref{eq:KF_joint2} reduce to 
\begin{align*}
    \bpy_{n+1} = \blG\mean_n, \quad \pCov_{n+1}^{\theta \by} = \pCov_{n+1} \blG^T,\quad  \textrm{and} \quad \pCov_{n+1}^{\by\by} = \blG  \pCov_{n+1} \blG^T + \Sigma_{\nu,n+1} \quad \textrm{where} \quad \blG = \begin{bmatrix}
      G  \\
      I \\
    \end{bmatrix}.
\end{align*}
We note that
\begin{equation}
    \label{eq:blG}
    \langle \blG v, \blG v \rangle \ge \|v\|^2.
\end{equation}
The update equations~\eqref{eq:KF_analysis_mean}\eqref{eq:KF_analysis_cov} become
\begin{subequations}
\begin{align}
        \mean_{n+1} &= \mean_{n} + \pCov_{n+1} \blG^T (\blG  \pCov_{n+1} \blG^T + \Sigma_{\nu,n+1})^{-1} (\by - \blG\mean_{n}), \label{eq:Lin_KF_analysis_mean}\\
         \Cov_{n+1}&= \pCov_{n+1} - \pCov_{n+1} \blG^T(\blG  \pCov_{n+1} \blG^T + \Sigma_{\nu,n+1})^{-1} \blG \pCov_{n+1}, \label{eq:Lin_KF_analysis_cov}
\end{align}
\end{subequations}
with $\pCov_{n+1} =  \Cov_{n} + \Sigma_{\omega,n+1}$.
We have the following theorem about the convergence of the algorithm:
\begin{theorem}
\label{th:lin_converge}
Assume that the error covariance matrices are as defined in~\cref{eq:hyperparameters} with $0 < \Dt < 1$ and that the prior covariance matrix $\Sigma_{0} \succ 0$ and initial covariance matrix $\Cov_{0} \succ 0.$
The iteration for the conditional mean $\mean_n$ and precision matrix $\Cov^{-1}_{n}$ characterizing the distribution of $\theta_n|Y_n$
converges exponentially fast to limit $\mean_{\infty}, \Cov^{-1}_{\infty}.$
Furthermore
the limiting mean $\mean_{\infty}$ 
and precision matrix $\Cov^{-1}_{\infty} = G^T\Sigma_{\eta}^{-1}G + \Sigma_0^{-1}$ are the 
posterior mean and precision matrix given by \eqref{eq:linear_posterior}. 
\end{theorem}

\begin{proof}
With the error covariance matrices defined in~\cref{eq:hyperparameters}, the update equation for $\{\Cov_n\}$ in \cref{eq:Lin_KF_analysis_cov} can be rewritten as 
\begin{equation}
\label{eq:Lin_KF_Cinv}
\begin{split}
    &\Cov_{n+1}^{-1} = \blG^T\Sigma_{\nu,n+1}^{-1}\blG + (\Cov_n + \Sigma_{\omega,n+1})^{-1}  = \dzh{\Dt} \Big(G^T\Sigma_{\eta}^{-1}G + \Sigma_0^{-1}\Big) + \dzh{(1 - \Dt)}\Cov_n^{-1}.\\
\end{split}
\end{equation}
We thus have a closed formula for $\Cov^{-1}_{n}$: 
\begin{equation}
\label{eq:beta_1}
    \Cov_n^{-1} = 
    \Big[1 - \dzh{(1 - \Dt)^n}\Big] \Big(G^T\Sigma_{\eta}^{-1}G + \Sigma_0^{-1}\Big) + \dzh{(1 - \Dt)^n} \Cov_0^{-1}.
\end{equation}
Since $0 < \Dt < 1$ this leads to the exponential convergence $\displaystyle \lim_{n\to \infty} \Cov_n^{-1} = G^T \Sigma_{\eta}^{-1} G + \Sigma_0^{-1} = \Cov_{\rm post}^{-1}$
given by \eqref{eq:linear_posterior}.

Since we have made a choice independent of $n$ we
write $\Sigma_{\nu} := \Sigma_{\nu,n+1}$. Thus
\cref{eq:Lin_KF_Cinv,eq:beta_1} lead to 
\begin{equation}
\label{eq:C_bound}
\begin{split}
 \blG^T\Sigma_{\nu}^{-1}\blG \preceq\Cov_{n+1}^{-1} \preceq \blG^T\Sigma_{\nu}^{-1}\blG + \Sigma_{+}
    \quad \textrm{where} \quad 
    \Sigma_{+} = 
    \dzh{\frac{1 - \Dt}{\Dt}}\blG^T\Sigma_{\nu}^{-1}\blG + \Cov_0^{-1}. 
\end{split}
\end{equation}
The update equation of $\mean_n$ in \cref{eq:Lin_KF_analysis_mean} can be rewritten as 
\begin{equation}
\label{eq:contracting-1}
\mean_{n+1} = \mean_{n} + \Cov_{n+1}\blG^T\Sigma_{\nu}^{-1}(\by - \blG\mean_{n}).
\end{equation}
Note that $B:=\blG^T\Sigma_{\nu}^{-1}\blG$ is symmetric and that, as a consequence of \eqref{eq:blG} together with the
fact that $\Sigma_{\nu} \succ 0$, it follows that $B \succ 0$;
thus we have that $I-\Cov_{n+1}B$ has the same
spectrum as $I-B^{\frac12}\Cov_{n+1}B^{\frac12}.$
Using the upper bound on $\Cov_{n+1}$ appearing in~\cref{eq:C_bound},
the spectral radius of the update matrix in \cref{eq:contracting-1} satisfies
\begin{equation}
\label{eq:cov-mean-par}
\begin{split}
\varrho(\I - \Cov_{n+1}\blG^T\Sigma_{\nu}^{-1}\blG) 
&= \varrho(\I - \Cov_{n+1}B) \\
&= \varrho\Big(\I - B^{\frac{1}{2}}\Cov_{n+1}B^{\frac{1}{2}}\Big) \\
&\leq
1 - \varrho\Big(B^{\frac{1}{2}}\big(B + \Sigma_{+}\big)^{-1}B^{\frac{1}{2}}\Big) \\
&= 1 -\epsilon_0,
\end{split}
\end{equation}
where $\epsilon_0 \in (0,1)$. Hence, we deduce that $\{\mean_{n}\}$ converges exponentially to the stationary point $\mean_{\infty}$, which satisfies $\displaystyle \blG^T\Sigma_{\nu}^{-1}(\by - \blG\mean_{\infty}) = 0$. Using the structure of $\blG$ and
$\Sigma_{\nu}$ the limiting mean can be written as the posterior
mean given in
\eqref{eq:linear_posterior}:  
\begin{align}
\label{eq:mean-inf}
    &m_{\infty} = r_0 + \Big(G^T \Sigma_{\eta}^{-1}  G + \Sigma_0^{-1}\Big)^{-1} G^T\Sigma_{\eta}^{-1}(y -  G r_0) = m_{\rm post}.
\end{align}
\end{proof}

\begin{remark}
Although this theorem applies only to the linear Gaussian setting
we note that the premise of matching only first and second order
moments is inherent to all Kalman methods. We demonstrate numerically in Section~\ref{sec:app} that application of the filtering methodology
based on the proposed choices of covariances leads to approximated mean and covariances which are accurate for nonlinear inverse problems.
\end{remark}

\begin{remark}
We note that the convergence of the means/covariances of the Kalman filter is 
a widely studied topic; 
and variants on some of our results can be obtained
from the existing literature, for example, the use of contraction mapping arguments to study
convergence of the Kalman filter is explored in~\cite{lancaster1995algebraic,bougerol1993kalman}.
\end{remark}

\subsection{Nonlinear Kalman Inversion Methodologies}
\label{sec:nonlinear_Kalman}
To make practical methods for solving nonlinear inverse problems~(\ref{eq:KI}) out of the foregoing, the expectations (integrals) appearing in the prediction step (\ref{eq:KF_pred})
as well as in the analysis step via \cref{eq:KF_joint2} need to be approximated appropriately. While
\cref{eq:KF_pred} can be implemented via a simple rescaling of the covariance matrix or ensemble, respectively, (we use both) the analysis step 
can be implemented using \emph{any} nonlinear Kalman filter 
(we use a variety). In the present work, we focus on both the unscented and ensemble Kalman filters, which lead to the Bayesian
implementations of unscented Kalman inversion (UKI), stochastic ensemble Kalman inversion (EKI), ensemble adjustment Kalman inversion (EAKI), and ensemble transform Kalman inversion (ETKI). We now detail these methods.
\footnote{Recall the discussion in Subsection \ref{ssec:1.3} about
distinction between optimization and Bayesian implementations of
all these methods.}

\subsubsection{Unscented Kalman Inversion~(UKI)}
\label{subsec:UKI}
UKI approximates the integrals in \cref{eq:KF_joint2} by means of deterministic quadrature rules; this is the idea of the unscented transform~\cite{julier1995new,wan2000unscented}. We now define
this precisely in the versions used in this paper.
\begin{definition}[Modified Unscented Transform~\cite{UKI}]
\label{def:unscented_transform}
Consider Gaussian random variable $\theta \sim \N(\mean, \Cov) \in \R^{N_{\theta}}$. Define $J$ sigma points $\{\theta^j\}_{j=0}^{J-1}$ according to the deterministic formulae
\begin{equation}
    \theta^0 = \mean \qquad \theta^j = \mean + [\sqrt{\Cov}] I_{N_{\theta}}[:,\, j] \quad (1\leq j\leq J - 1);
\end{equation}
here $[\sqrt{\Cov}]$ is the Cholesky factor of $\Cov$ and 
$I_{N_{\theta}}[:, j]$ is the $j^{th}$ column of the matrix $I_{N_{\theta}} \in \R^{N_\theta \times (J-1)}$.
Consider any two real vector-valued functions 
$\bG_1(\cdot)$ and $\bG_2(\cdot)$ acting on $\R^{N_\theta}.$ Using the sigma points we may define a quadrature rule approximating
the mean and covariance of the  random variables $\bG_1(\theta)$ and $\bG_2(\theta)$ as follows:  
\begin{equation}
\label{eq:ukf-original}
\begin{split}
    \E[\bG_i(\theta)] \approx \bG_i(\theta^0)\qquad 
    \textrm{Cov}[\bG_1(\theta), \bG_2(\theta]  \approx \sum_{j=1}^{J-1} a \Big(\bG_1(\theta^j) - \E\bG_1(\theta)\Big)\Big(\bG_2(\theta^j) - \E\bG_2(\theta)\Big)^T. 
\end{split}
\end{equation}

In the present work, we consider the following two variants,
\begin{itemize}
    \item UKI-1 $(J = N_{\theta}+2)$~\cite{julier2002reduced,moireau2011reduced}. $I_{N_{\theta}}$ is defined recursively as
\begin{align}
& I_1 = 
    \begin{bmatrix}
    -\frac{1}{\sqrt{2 a}} & \frac{1}{\sqrt{2 a}}
\end{bmatrix} \\
& I_d = 
    \begin{bmatrix}
     & & & 0 \\
     & I_{d-1}& & \vdots\\
     & & & 0\\
    \frac{1}{\sqrt{ad(d+1)}} & \dots & \frac{1}{\sqrt{ad(d+1)}} & \frac{-d}{\sqrt{ad(d+1)}}
    \end{bmatrix}, \qquad 2\leq d \leq N_{\theta}
\end{align}
and the weight parameter is chosen as $a = \frac{N_{\theta}}{4(N_{\theta}+1)}$. 
    \item UKI-2 $(J = 2N_{\theta}+1)$~\cite{UKI}. $I_{N_{\theta}}$ is defined as
\begin{align}
& I_{N_{\theta}} = 
    \begin{bmatrix}
     \frac{1}{\sqrt{2a}}  & & & & -\frac{1}{\sqrt{2a}}          & & &       \\
     &  \frac{1}{\sqrt{2a}}& & & &  -\frac{1}{\sqrt{2a}}& & &\\
     & & \ddots&   &     &    &  \ddots&  &     \\
     &  &    &  \frac{1}{\sqrt{2a}}  &  &  &       & -\frac{1}{\sqrt{2a}} &
    \end{bmatrix}
\end{align}
and the weight parameter is chosen as $a = \max \{\frac{1}{8}, \frac{1}{2N_{\theta}} \}$. 
    
\end{itemize}

\end{definition}

Consider the Gaussian approximation algorithm defined by \cref{eq:KF_pred,eq:KF_joint,eq:KF_joint2,eq:KF_analysis_mean,eq:KF_analysis_cov}.
By utilizing the aforementioned quadrature rule, the iteration procedure of the UKI becomes:
\begin{itemize}
\item Prediction step : 
    \begin{align}
        &\pmean_{n+1} = \mean_{n} \qquad \pCov_{n+1} = \dzh{\frac{1}{1 - \Dt}}\Cov_n.
    \end{align}
    \item Generate sigma points :
    \begin{align}
    \pp_{n+1}^0 = \pmean_{n+1},
    \qquad \pp_{n+1}^j = \pmean_{n+1} + [\sqrt{\pCov_{n+1}}] I_{N_{\theta}}[:,\, j] \quad (1\leq j\leq J-1).
    \end{align}
\item Analysis step :
  \begin{subequations}
  \label{eq:UKI-sample}
  \begin{align}
        &\bpy^j_{n+1} = \bG(\pp^j_{n+1}) \quad (0\leq j \leq J-1), \label{eq:UKI-forward}\\
        &\bpy_{n+1} = \bpy^0_{n+1},\nonumber\\
         &\pCov^{\theta \by}_{n+1} = \sum_{j=1}^{J-1}a
        (\pp^j_{n+1} - \pmean_{n+1} )(\bpy^j_{n+1} - \bpy_{n+1})^T, \nonumber\\
        &\pCov^{\by\by}_{n+1} = \sum_{j=1}^{J-1}a
        (\bpy^j_{n+1} - \bpy_{n+1} )(\bpy^j_{n+1} - \bpy_{n+1})^T + \Sigma_{\nu,n+1},\nonumber\\
        &\mean_{n+1} = \pmean_{n+1} + \pCov^{\theta \by}_{n+1}(\pCov^{\by\by}_{n+1})^{-1}(\by - \bpy_{n+1}),\nonumber\\
        &\Cov_{n+1} = \pCov_{n+1} - \pCov^{\theta \by}_{n+1}(\pCov^{\by\by}_{n+1})^{-1}{\pCov^{\theta \by}_{n+1}}{}^{T}.\nonumber
    \end{align}
    \end{subequations}
\end{itemize}

\subsubsection{Ensemble Kalman Inversion}
\label{sec:app:ensemble_method}
Ensemble Kalman inversion represents the distribution at each iteration by an ensemble of parameter estimates $\{\theta_n^j\}_{j=1}^{J}$and approximates the integrals in \cref{eq:KF_joint2} empirically. We describe three variants on this idea.

\paragraph{{\textbf{Stochastic Ensemble Kalman Inversion~(EKI)}}}
The perturbed observations form of the ensemble Kalman filter
\cite{evensen2009data} is applied to the extended mean-field dynamical system \eqref{eq:evolve}-\eqref{eq:obser}, which leads to the EKI:
\begin{itemize}
  \item Prediction step : 
  \begin{align}
  \label{eq:eki_pred}
  \pmean_{n+1} = \mean_n,\qquad
  \pp_{n+1}^{j} = \pmean_{n+1} + \dzh{\sqrt{\frac{1}{1 - \Dt}}} (\theta_{n}^{j} - \mean_n).
 \end{align}
  \item  Analysis step : 
  \begin{subequations}
    \label{eq:EKI-analysis}
  \begin{align}
         &\bpy_{n+1}^{j} = \bG(\pp_{n+1}^{j})  \qquad \bpy_{n+1} = \frac{1}{J}\sum_{j=1}^{J}\bpy_{n+1}^{j},\label{eq:EKI-Fx-analysis}\\
        &\pCov_{n+1}^{\theta \by} = \frac{1}{J-1}\sum_{j=1}^{J}(\pp_{n+1}^{j} - \pmean_{n+1})(\bpy_{n+1}^{j} - \bpy_{n+1})^T,  \\
        &\pCov_{n+1}^{\by\by} = \frac{1}{J-1}\sum_{j=1}^{J}(\bpy_{n+1}^{j} - \bpy_{n+1})(\bpy_{n+1}^{j} - \bpy_{n+1})^T +\Sigma_{\nu,n+1}, \label{eq:EKI-Cxx-analysis}\\
        &\theta_{n+1}^{j} = \pp_{n+1}^{j} + \pCov_{n+1}^{\theta \by}\left(\pCov_{n+1}^{\by\by}\right)^{-1}(\by - \bpy_{n+1}^{j} - \nu_{n+1}^{j}), \label{eq:EKI-theta-analysis}\\
        &\mean_{n+1} = \frac{1}{J} \sum_{j=1}^{J} \theta_{n+1}^{j}.
\end{align}
\end{subequations}
\end{itemize}
Here the superscript $j = 1,\cdots,\ J$ is the ensemble particle index, and  $\nu_{n+1}^{j} \sim \N(0, \Sigma_{\nu,n+1})$ are independent and identically distributed random variables. The prediction step ensures the exactness of the predictive covariance~\cref{eq:KF_pred}.

\begin{remark}
The prediction step~\eqref{eq:eki_pred} is inspired by square root Kalman filters~\cite{anderson2001ensemble,bishop2001adaptive,tippett2003ensemble, wang2003comparison} and covariance inflation~\cite{anderson2007adaptive}; these
methods are designed to ensure that the mean and covariance of $\{\pp_{n+1}^{j}\}$ match $\pmean_{n+1}$ and $\pCov_{n+1}$ exactly. This is different from traditional stochastic ensemble Kalman inversion implementation, where i.i.d. Gaussian noises $\omega^j_{n+1} \sim \N(0, \Sigma_{\omega, n+1})$ are added. 
In the analysis step~\eqref{eq:EKI-analysis}, we add noise in the $\{\theta_{n+1}^j\}$ update~\eqref{eq:EKI-theta-analysis} instead of the $\{\bpy_{n+1}^{j}\}$ evaluation~\eqref{eq:EKI-Fx-analysis}; this ensures that $\pCov_{n+1}^{\by\by}$~\eqref{eq:EKI-Cxx-analysis} is symmetric positive definite. 
\end{remark}

\begin{remark}
As a precursor to understanding the adjustment and transform
filters which follow this subsection, we show that the EKI
 does not exactly replicate the covariance update equation~\eqref{eq:KF_analysis_cov}. To this end,
denote  the matrix square roots $\pCovZ_{n+1},\, \CovZ_{n+1} \in \R^{N_{\theta}\times J}$ of $\pCov_{n+1},\,\Cov_{n+1}$ and $\pGCovZ_{n+1}$ as follows:
\begin{equation}
\label{eq:EKF_square_root}
\begin{split}
    \pCovZ_{n+1} &= \frac{1}{\sqrt{J-1}}\Big(\pp_{n+1}^{1} - \pmean_{n+1}\quad \pp_{n+1}^{2} - \pmean_{n+1}\quad...\quad\pp_{n+1}^{J} - \pmean_{n+1} \Big),\\
     \CovZ_{n+1} &= \frac{1}{\sqrt{J-1}}\Big(\theta_{n+1}^{1} - \mean_{n+1}\quad \theta_{n+1}^{2} - \mean_{n+1}\quad...\quad\theta_{n+1}^{J} - \mean_{n+1} \Big),\\
     \pGCovZ_{n+1} &= \frac{1}{\sqrt{J-1}}\Big(\bpy_{n+1}^{1} - \bpy_{n+1}\quad \bpy_{n+1}^{2} - \bpy_{n+1}\quad...\quad\bpy_{n+1}^{J} - \bpy_{n+1} \Big).
     \end{split}
\end{equation}
Then the covariance update equation~\eqref{eq:KF_analysis_cov} does not hold exactly:
\begin{equation}
\begin{aligned}
        \pCov_{n+1} - \pCov_{n+1}^{\theta \by}(\pCov_{n+1}^{\by\by})^{-1} {\pCov_{n+1}^{\theta \by}}{}^{T} &= \pCovZ_{n+1}\pCovZ_{n+1}^T - \pCovZ_{n+1}\pGCovZ_{n+1}^T (\pGCovZ_{n+1}\pGCovZ_{n+1}^T + \Sigma_{\nu,n+1})^{-1}\pGCovZ_{n+1}\pCovZ_{n+1}^T \\
        &\neq Z_{n+1} Z_{n+1}^T = C_{n+1}.
\end{aligned}
\end{equation}
\end{remark}

\paragraph{\textbf{Ensemble Adjustment Kalman Inversion~(EAKI)}}
Following the ensemble adjustment Kalman filter proposed in~\cite{anderson2001ensemble}, the analysis step updates particles deterministically with a pre-multiplier $A$,
\begin{equation}
\theta_{n+1}^j - \mean_{n+1} = A (\pp_{n+1}^{j} - \pmean_{n+1}). \\
\end{equation}
Here $A = P \pD^{\frac{1}{2}} U D^{\frac{1}{2}}\pD^{-\frac{1}{2}}P^T $ with 
\begin{equation}
\label{eq:eaki-svd}
\begin{split}
   \textrm{SVD :}     \quad       &\pCovZ_{n+1} =  P \pD^{\frac{1}{2}} V^T,\\
   \textrm{SVD :}     \quad      &V^T\Big(\I + \pGCovZ_{n+1}^T \Sigma_{\nu,n+1}^{-1}  \pGCovZ_{n+1}\Big)^{-1} V = U D U^T,
\end{split}
\end{equation}
where both $\pD$ and $D$ are non-singular diagonal matrices, with dimensionality rank($\pCovZ_{n+1}$), and $\pCovZ_{n+1}$ and $\pGCovZ_{n+1}$ are defined in \cref{eq:EKF_square_root}.
The analysis step becomes:
\begin{itemize}
  \item  Analysis step : 
\begin{subequations}
\begin{align}
        &\mean_{n+1} = \pmean_{n+1} + \pCov_{n+1}^{\theta \by}\left(\pCov_{n+1}^{\by\by}\right)^{-1}(\by - \bpy_{n+1}),\label{eq:EAKI-mean-analysis}\\
        &\theta_{n+1}^{j} = \mean_{n+1} + A(\pp_{n+1}^{j} - \pmean_{n+1}).\label{eq:EAKI-particle-analysis}
\end{align}
\end{subequations}
\end{itemize}

\begin{remark}
It can be verified that the covariance update equation~\eqref{eq:KF_analysis_cov} holds:
\begin{equation}
\label{eq:cov-eaki-2}
    \begin{split}
        \Cov_{n+1} &= \CovZ_{n+1} \CovZ_{n+1}^T \\
                   &= A\pCovZ_{n+1}  \pCovZ_{n+1}^T A^T \\
                   &= P \pD^{\frac{1}{2}}  U D U^{T} \pD^{\frac{1}{2}}P \\
                   &= \pCovZ_{n+1}\Big(\I +  \pGCovZ_{n+1}^T  \Sigma_{\nu,n+1}^{-1}\pGCovZ_{n+1}\Big)^{-1}\pCovZ_{n+1}^T \\
                   &= \pCovZ_{n+1}\Big(\I -  \pGCovZ_{n+1}^T (\pGCovZ_{n+1}\pGCovZ_{n+1}^T + \Sigma_{\nu,n+1})^{-1}\pGCovZ_{n+1}\Big)\pCovZ_{n+1}^T \\
                   &= \pCov_{n+1} - \pCov_{n+1}^{\theta \by}(\pCov_{n+1}^{\by\by})^{-1} {\pCov_{n+1}^{\theta \by}}{}^{T}.
    \end{split}
\end{equation}
\end{remark}

\paragraph{\textbf{Ensemble Transform Kalman Inversion~(ETKI)}}

Following the ensemble transform Kalman filter proposed in~\cite{bishop2001adaptive,tippett2003ensemble, wang2003comparison}, the analysis step updates particles deterministically with a post-multiplier $T$,
\begin{equation}
\CovZ_{n+1} = \pCovZ_{n+1} T.\\
\end{equation}
Here $T = P(\Gamma + I)^{-\frac{1}{2}}P^T$, with 
\begin{equation}
\textrm{SVD:} \quad \pGCovZ_{n+1} \Sigma_{\nu,n+1}^{-1} \pGCovZ_{n+1} = P\Gamma P^T.
\end{equation}
The analysis step becomes:
\begin{itemize}
  \item  Analysis step : 
\begin{subequations}
\begin{align}
        &\mean_{n+1} = \pmean_{n+1} + \pCov_{n+1}^{\theta \by}\left(\pCov_{n+1}^{\by\by}\right)^{-1}(\by - \bpy_{n+1}),\label{eq:ETKI-mean-analysis}\\
        &\CovZ_{n+1} = \pCovZ_{n+1} T. \label{eq:ETKI-particle-analysis}
\end{align}
\end{subequations}
\end{itemize}

\begin{remark}
It can be verified that the covariance update equation~\eqref{eq:KF_analysis_cov} holds:
\begin{equation}
\label{eq:cov-etki-2}
    \begin{split}
        \Cov_{n+1} &= \CovZ_{n+1} \CovZ_{n+1}^T  \\
                   &= \pCovZ_{n+1} T T^T \pCovZ_{n+1}^T \\
                   &= \pCovZ_{n+1}\Big(\I +  P\Gamma P^T\Big)^{-1}\pCovZ_{n+1}^T \\
                   &= \pCovZ_{n+1}\Big(\I +  \pGCovZ_{n+1}^T  \Sigma_{\nu,n+1}^{-1}\pGCovZ_{n+1}\Big)^{-1}\pCovZ_{n+1}^T \\
                   &= \pCovZ_{n+1}\Big(\I -  \pGCovZ_{n+1}^T (\pGCovZ_{n+1}\pGCovZ_{n+1}^T + \Sigma_{\nu,n+1})^{-1}\pGCovZ_{n+1}\Big)\pCovZ_{n+1}^T \\
                   &= \pCov_{n+1} - \pCov_{n+1}^{\theta \by}(\pCov_{n+1}^{\by\by})^{-1} {\pCov_{n+1}^{\theta \by}}{}^{T}.
    \end{split}
\end{equation}
\end{remark}

Particles~(ensemble members) updated by the basic form of the EKI algorithm through iterates are confined to
the linear span of the initial ensemble \cite{emerick2013investigation,iglesias2013ensemble}. 
The same is true for both EAKI and ETKI: 

\begin{lemma}
\label{lemma}
For both EAKI and ETKI, all particles lie in the linear space $\mathcal{A}$ spanned by $m_0$ and the column vectors of $Z_0$.  
\end{lemma}
\begin{proof}
We will prove that $m_n$ and column vectors of $\CovZ_n$ are in $\mathcal{A}$ by induction. We assume this holds for all $n \leq k$.
Since $\pmean_{k+1} = \mean_{k}$ and $\pCovZ_{k+1} = \dzh{\frac{1}{1 - \Dt}}\CovZ_{k}$ (See \cref{eq:eki_pred}), $\pmean_{k+1}$ and column vectors of $\pCovZ_{k+1}$ are in $\mathcal{A}$.
Combining the mean update equations~(\ref{eq:EAKI-mean-analysis},\ref{eq:ETKI-mean-analysis}) and the fact that $\pCov_{k+1}^{\theta \by} = \pCovZ_{k+1} \pGCovZ_{k+1}^T$, we have $\mean_{k+1}$ is in $\mathcal{A}$.
For EAKI, since the pre-mulitiplier $A = P \pD^{\frac{1}{2}} U D^{\frac{1}{2}} \pD^{-\frac{1}{2}}P^T$, and $P$ is the left compact singular matrix of $\pCovZ_{k+1}$, it follows that the column vectors of $A$ lie in $\mathcal{A}$; furthermore, the square root matrix update equation~(\ref{eq:EAKI-particle-analysis}), $\CovZ_{k+1} = A\pCovZ_{k}$, has
implication that the column vectors of $\CovZ_{k+1}$ lie in $\mathcal{A}$.
For the ETKI, the square root matrix update equation~(\ref{eq:ETKI-particle-analysis}) implies  that the column vectors of $\CovZ_{k+1}$ lie in $\mathcal{A}$.
Since $m_n$ and column vectors of $\CovZ_n$ are in $\mathcal{A}$, so are the particles $\{\theta_n^{j}\}_{j=1}^{J}$.
\end{proof}

\section{Variants on the Basic Algorithm}
\label{sec:Variance}

In this section, we introduce three strategies
to make the novel mean-field based methodology more efficient, robust and
widely applicable in real large-scale problems. In Subsection
\ref{subs:low-rank} we introduce low-rank approximation, in which the parameter space is restricted to a low-rank space induced by the prior;
Subsection \ref{subs:bi-fidelity} introduces a bi-fidelity approach
in which multifidelity models are used for different ensemble members;
and box constraints to enforce pointwise bounds on $\theta$ are
introduced in Subsection \ref{subs:box}.

\subsection{Low-Rank Approximation}
\label{subs:low-rank}

When using ensemble methods for state estimation, the dimension
of the ensemble space needed for successful state estimation may
be much smaller than $N_\theta;$ a useful rule of thumb is that
the ensemble space needs to be rich enough to learn about
the unstable directions in the system. When using ensemble methods
for inversion the situation is not so readily understood.
The EKI algorithm  presented here is limited to finding solutions 
in the linear span of the initial ensemble \cite{emerick2013investigation,iglesias2013ensemble} and we have
highlighted a similar property for the EAKI and ETKI in Lemma \ref{lemma}.
Whilst localization is often used to break this property \cite{chen2017localization} its
use for this purpose is somewhat ad hoc. In this work we do not
seek to break the subspace property. Indeed here we exploit
low rank approximation within ensemble inversion techniques, a methodology which leads to
solutions restricted to the linear span of a small number of
dominant modes defined by the prior distribution.

\Cref{th:lin_converge} requires that the initial covariance matrix $\Cov_{0} \succ 0$ be strictly
positive definite. 
To satisfy the assumption, the UKI requires $N_{\theta}+2$ or $2N_{\theta}+1$ forward problem evaluations and the storage of an $N_{\theta} \times N_{\theta}$ covariance matrix, and the EKI, EAKI and ETKI require $\bigO(N_{\theta})$ forward problem evaluations and the storage of $\bigO(N_{\theta})$ parameter estimates; some of
the implications of these effects are numerically verified in \cref{sec:app-Hilbert}. Therefore, they are unaffordable for field inversion problems, where $N_{\theta}$ is large, typically from discretization of the $N_{\theta} = \infty$ limit.
However, many physical phenomena or systems exhibit large-scale structure or finite-dimensional attractors, and in such situations the model error covariance matrices are generally low-rank; these low-rank spaces are spanned by, for example, the dominant Karhunen-Lo\`{e}ve modes for random fields~\cite{ghanem2003stochastic,stuart2010inverse} or the dominant spherical harmonics space on the sphere~\cite{anderson2001ensemble,ehrendorfer2011spectral}.
We introduce a reparameterization strategy for this framework in order to leverage such low-rank structure when present, and thereby to reduce both computational and storage costs.

Given the prior distribution $\N(r_0, \Sigma_0)$, we assume $\Sigma_0$ is a low-rank matrix with the truncated singular value decomposition
$$\Sigma_0 \approx U D_0 U^T.$$
Here $U = \{u_1, u_2,\cdots, u_{N_r}\}$ is the $N_r$-dominant singular vector matrix and $D_0$ is the singular value matrix.
The discrepancy $\theta - r_0$ is assumed to be well-approximated in the linear space spanned by column vectors of $U$. Hence, the unknown parameters can be reparameterized as follows: 
\begin{equation*}
    \theta = r_0 + \sum_{i=1}^{N_r}\tau_{(i)} u_i.
\end{equation*}
The aforementioned algorithm is then applied to solve for the vector $\tau = [\tau_{(1)}, \tau_{(2)}, \cdots, \tau_{(N_r)}]^T$, which has prior mean $0$ and prior covariance $D_0$.
This reduces the computation and memory cost from $\bigO(N_{\theta})$ and $\bigO(N_\theta^2)$ to $\bigO(N_r)$ and $\bigO(N_rN_{\theta})$, where $N_{r}$ is the rank of the covariance matrix.

More advanced approaches to extracting the low-rank space exist,  including active subspace methods~\cite{constantine2014active} and likelihood-informed subspace methods~\cite{cui2014likelihood,cui2021unified}; however, they all require derivatives and so we do not pursue them here.

\subsection{Bi-Fidelity Approach}
\label{subs:bi-fidelity}
For large-scale scientific or engineering problems, even with a small parameter number $N_\theta$~(or rank number $N_r$), the computational cost associated with these $\bigO(N_{\theta})$~(or $\bigO(N_{r})$) forward model evaluations can be intractable; for example the number of parameters may be small, but the parameter-to-data map may require
evaluation of a large, complex model. The bi-fidelity or multilevel strategy~\cite{heinrich2001multilevel,giles2015multilevel,gao2021bi,fairbanks2020bi} is widely used to accelerate sampling-based methods; in
particular it has been introduced in the context of ensemble
methods in \cite{hoel2016multilevel} and see \cite{chada2020multilevel} for
a recent overview of developments in this direction.

We employ a particular bi-fidelity approach for the UKI algorithm. In this
approach, low-fidelity models can be used to speed up forward model evaluations as follows. Consider \cref{eq:UKI-forward};
evaluation of the mean $\G(\hat{\theta}_{n+1}^0)$ 
can be performed using a high-fidelity model; meanwhile the other $J-1$ forward evaluations employed for covariance estimation, $\{\G(\hat{\theta}_{n+1}^j)\}_{j=1}^{J-1}$, can use low-fidelity models.

\subsection{Box Constraints}
\label{subs:box}

Adding constraints to the parameters~(for example, dissipation is non-negative) significantly improves the robustness of Kalman inversion~\cite{albers2019ensemble,chada2019incorporation,schneider2020ensemble}.
In this paper, there are occasions where
we impose element-wise box constraints of the form
\begin{equation*}
    0 \leq \theta  \quad \textrm{ or }\quad \theta_{min} \leq \theta \leq \theta_{max}.
\end{equation*}
These are enforced by change of variables writing $\theta = \varphi(\tilde{\theta})$
where, for example, respectively,
\begin{equation*}
   \varphi(\tilde{\theta}) = \exp(\tilde{\theta}) \quad \textrm{ or } \varphi(\tilde{\theta}) = \theta_{min} +  \frac{\theta_{max}-\theta_{min}}{1 + \exp(\tilde\theta)}.
\end{equation*}
The inverse problem is then reformulated as 
\begin{equation*}
    y = \G(\varphi(\tilde\theta)) + \eta.
\end{equation*}
and the proposed Kalman inversion methods are employed with $\G \mapsto \G \circ \varphi.$

\section{Numerical Experiments}
\label{sec:app}

In this section, we present numerical experiments demonstrating application of filtering methods to the novel mean-field dynamical system (\cref{eq:evolve,eq:obser,eq:hyperparameters}) introduced in this paper
\footnote{\dzh{We fix $\Dt = 1/2$ based on the parameter study presented in \ref{appendix:CTL} and iterate $\bigO(10)$ iterations to demonstrate convergence. In practice, adaptive time stepping and increment-based stopping criteria can be applied.}}; the goal is to approximate the posterior distribution of  unknown parameters or fields. The first subsection lists the five test problems, and the subsequent subsections consider them each in turn. In summary, our findings are as follows:
\footnote{The footnote from Subsection \ref{ssec:1.4}, appearing
before the bulleted list of contributions, applies here too.}

\begin{itemize}
    \item The proposed Kalman inversion methods based on  (\cref{eq:evolve,eq:obser,eq:hyperparameters})
    are more efficient than transport/coupling methods
    based on \eqref{eq:KL_dynamics}~(i.e., iterative Kalman filter methods) 
    on all the examples we consider. They remove the sensitivity to the initialization and,
    relatedly, they converge exponentially fast.
    \item The proposed Kalman inversion methods with deterministic treatment of stochastic terms, specifically UKI and EAKI, outperform other methods with stochastic treatments, such as EKI, EKS (with Euler-Maruyama discretization) and CBS. They do not suffer from the presence of noisy fluctuations and achieve convergence for both linear and nonlinear problems.
    \item The methodology is implementable for large-scale parameter identification problems, such as those arising in climate models.
    \end{itemize}

\subsection{Overview of Test Problems}
\label{ssec:Test Problems}

The five test problems considered are:

\begin{enumerate}
    \item Linear-Gaussian $2$-parameter model problem: this problem serves as a proof-of-concept example, which demonstrates the convergence of the mean and the covariance as analyzed in Subsection~\ref{ssec:lin}.
    \item Nonlinear 2-parameter two-point boundary value problem: this example appears in
    \cite{ernst2015analysis} an important paper which demonstrates that the mean field
    limit of ensemble Kalman inversion methods may be far from the true
    posterior; it is also used as a test case in several other papers, such as~\cite{garbuno2020interacting,herty2018kinetic}. We show that by applying
    Kalman filtering techniques to the extended mean-field dynamical system (\cref{eq:evolve,eq:obser,eq:hyperparameters}),
    we obtain methods which obtain accurate posterior approximation on this problem.
     \item Hilbert matrix problem: this high dimensional linear-Gaussian problem demonstrates the ability of the proposed Kalman inversion
     methodology to solve ill-conditioned inverse problems.
    In addition to testing the novel mean-field approach
    introduced in this paper, we also 
    study the effect of the ensemble size on ensemble Kalman inversion, and in particular, the ensemble adjustment/transform Kalman inversions are examined
    in this context.
    \item Darcy flow inverse problem: this is an infinite dimensional field inversion problem (see \cite{dashti2013bayesian} and the references therein); in addition to testing the novel mean-field approach
    introduced in this paper, we also demonstrate the low-rank approximation strategy in Subsection~\ref{subs:low-rank}.
    \item Idealized global climate model: this 3D Navier-Stokes problem,
    see \cite{UKI} for background and references;
    in addition to testing the novel mean-field approach
    introduced in this paper, we also demonstrate the bi-fidelity approach introduced in Subsection~\ref{subs:bi-fidelity}.
\end{enumerate}

For the first and third problems, the Gaussian structure means
that they are exactly solvable and this provides a benchmark against
which we compare various methods.
Markov Chain Monte Carlo methods~(MCMC), specifically the random walk Metropolis~\cite{gelman1997weak} and preconditioned Crank–Nicolson~\cite{cotter2013mcmc} methods, 
are used as the benchmark for the second and fourth problems respectively.
Problem five is too large for the use of MCMC, and showcases
the potential of the methodology studied here to solve problems
otherwise beyond reach.

In the first two tests, we compare the proposed Kalman inversion methods 
(EKI, UKI, EAKI, ETKI  applied to \cref{eq:evolve,eq:obser,eq:hyperparameters} with other recently proposed Gaussian approximation algorithms, including 
the ensemble Kalman sampler~(EKS)~\cite{garbuno2020interacting,nusken2019note}
\footnote{We follow the implementation in~\cite{garbuno2020interacting}, which employs adaptive time-stepping.},
and the  consensus-based sampler~(CBS)~\cite{pinnau2017consensus, carrillo2018analytical,carrillo2021consensus_constraints,carrillo2021consensus_ML,carrillo2021consensus_sampling}\footnote{We follow the implementation in~\cite{carrillo2021consensus_sampling}, setting $\alpha=0.0$ and
adaptively updating $\beta$  with $\mu = 0.5$.}. We also compare with
variants of iterative Kalman filter methods, which seek to deform
the prior into the posterior in one time unit (transport/coupling) 
using a finite number of intermediate steps (see~\ref{appendix:CTEKI})
based on \eqref{eq:KL_dynamics}. They include iterative unscented Kalman filters (IUKF-1, IUKF-2), iterative ensemble Kalman filter (IEnKF), iterative ensemble adjustment Kalman filter~(IEAKF), and iterative ensemble transform Kalman filter~(IETKF)~\cite{bocquet,oliver2008inverse,emerick2013investigation,chen2012ensemble}\cite[Algorithm 3]{ding2020ensemble}. 
Having shown the superiority of filtering based on our novel mean field
dynamical system, we consider only this approach in the remaining examples.
In the third test, we study the effect of the ensemble size on the proposed Kalman
inversion methods, in particular for the EAKI/ETKI approaches to filtering, comparing with the UKI.
In the fourth and fifth tests, we demonstrate the effectiveness of the proposed Kalman inversion methods for large-scale inverse problems and incorporate the low-rank and bi-fidelity approaches.

\subsection{Linear 2-Parameter Model Problem}
\label{ssec:app:add}
We consider the $2$-parameter linear inverse problem~\cite{UKI}
of finding $\theta \in \R^2$ from $y$ given by
\begin{align}
    y = G\theta + \eta.
\end{align}
Here the observation error noise
is $\eta \sim \N(0, 0.1^2\I)$. We explore the following two scenarios
\begin{itemize}
\item over-determined system
    \begin{align}
    y = \begin{bmatrix}
        3 \\
        7 \\
        10
       \end{bmatrix}~~~
       G = \begin{bmatrix}
        1 & 2\\
        3 & 4 \\
        5 & 6
    \end{bmatrix}~~~
    \rho_{\rm prior} \sim \N(0, \I);
   \end{align} 
    \item under-determined system
    \begin{align}
    y = \begin{bmatrix}
        3
       \end{bmatrix}~~~
       G = \begin{bmatrix}
        1 & 2\\
        \end{bmatrix}~~~
    \rho_{\rm prior} \sim \N(0, \I).
   \end{align} 
\end{itemize}

\begin{figure}[ht]
\centering
    \includegraphics[width=0.8\textwidth]{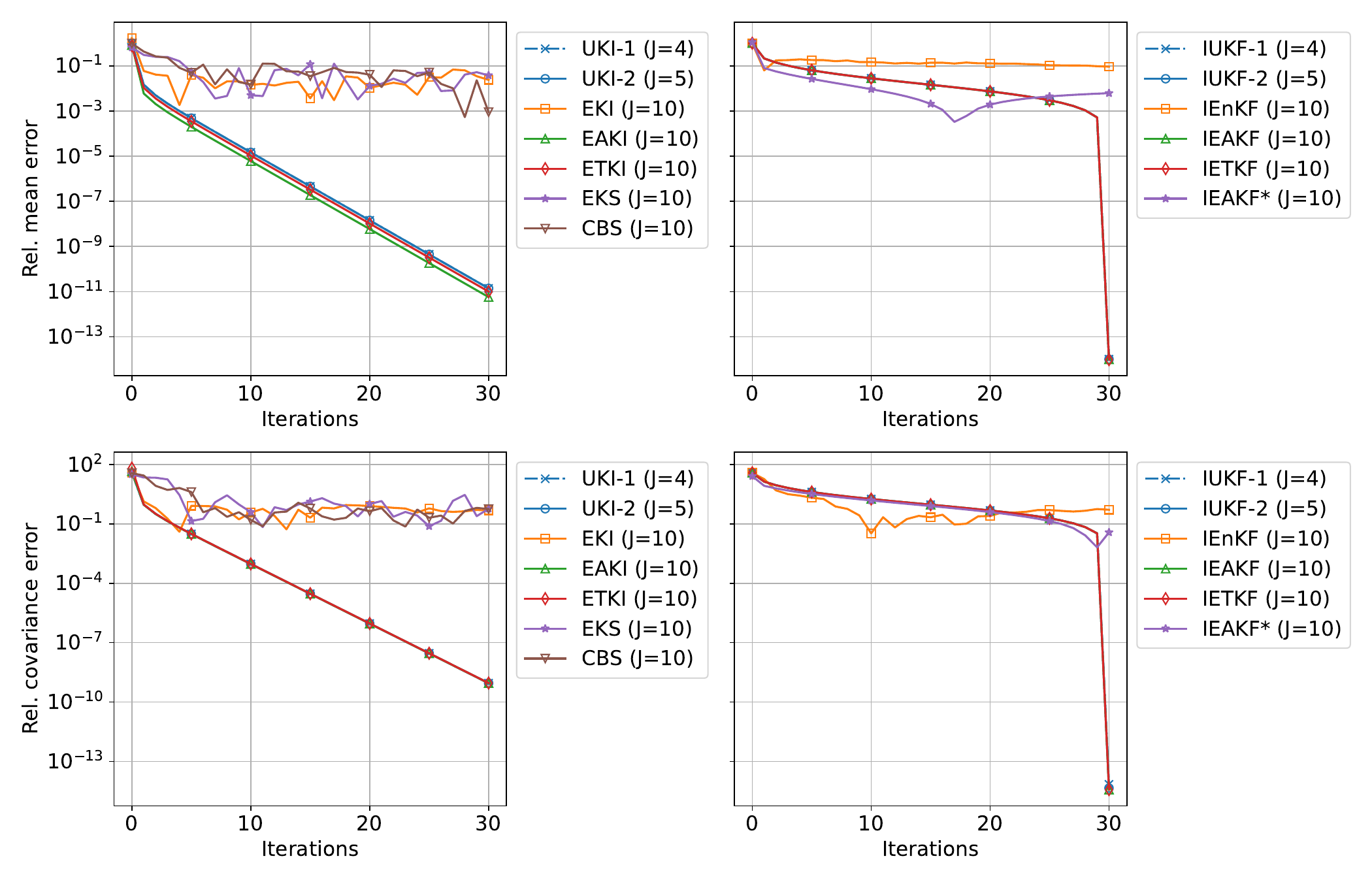}
    \caption{Linear 2-parameter model problems: convergence of posterior mean (top) and posterior covariance (bottom) for the over-determined system.}
    \label{fig:linear-well}
\end{figure}

\begin{figure}[ht]
\centering
    \includegraphics[width=0.8\textwidth]{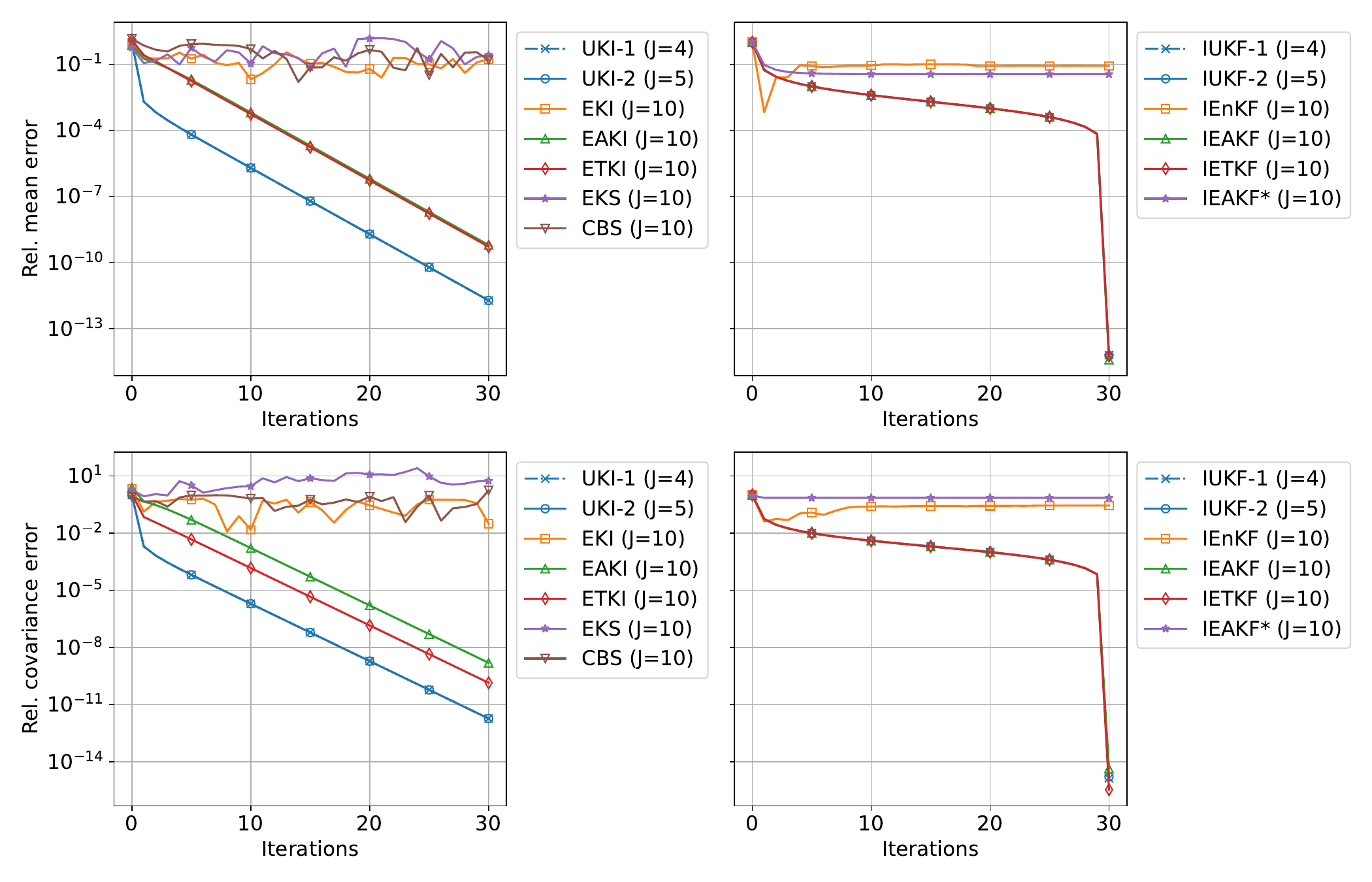}
    \caption{Linear 2-parameter model problems: convergence of posterior mean (top) and posterior covariance (bottom) for the under-determined system (bottom).}
    \label{fig:linear-under}
\end{figure}

We apply various Kalman inversions to our
proposed novel mean-field dynamical system, including UKI-1~($J=4$), UKI-2~($J=5$), EKI, EAKI,  and ETKI; we compare with pre-existing coupling/transport based
iterative Kalman filters, including IUKF-1~($J=4$), IUKF-2~($J=5$), IEnKF, IEAKF and IETKF all with $J = 10$ ensemble members; and we compare with  EKS, and CBS,  again all with $J = 10$.
All algorithms are initialized at the prior distribution; note however that
the methods we introduce in this paper, and EKS and CBS, do not require this and indeed are robust to the use of different initializations.
The iterative Kalman filters are discretized with $\Dt = \frac{1}{30}$, and further correction~(See \ref{appendix:Gaussian_Init}) is applied on the initial ensemble members for the exactness of the initialization, except for the IEAKF$^{*}$.
Since the posterior distribution is Gaussian, we can compute the reference distribution analytically. The convergence of the posterior mean and posterior covariance are reported in~\cref{fig:linear-well,fig:linear-under}. 
Because we use the same number of steps for all algorithms, and commensurate numbers of particles, the evaluation cost of all the methods
studied are comparable; the size of the error discriminates between
them.

For both scenarios, UKI-1, UKI-2, EAKI, and ETKI converge exponentially fast.
IUKF-1, IUKF-2, IEAKF, and IETKF reach exact posterior mean and covariance matrix at $T=1$.
However, IEAKF$^{*}$ does not converge due to the error introduced in the initialization.
EKI and IEnKF do not converge, and suffer from the presence of
random noise introduced in the analysis step. 
EKS and CBS do not converge, and suffer from the presence of random noise and the finite ensemble size.

\subsection{Nonlinear 2-Parameter Model Problem}
\label{ssec:app:2-p}

Consider the one-dimensional elliptic boundary-value problem
\begin{equation}
    -\frac{d}{dx}\Big(\exp(\theta_{(1)}) \frac{d}{dx}p(x)\Big) = 1, \qquad x\in[0,1]
\end{equation}
with boundary conditions $p(0) = 0$ and $p(1) = \theta_{(2)}$. The solution for this problem is given by  
\begin{equation}
    p(x) = \theta_{(2)} x + \exp(-\theta_{(1)})\Big(-\frac{x^2}{2} + \frac{x}{2}\Big).
\end{equation}
The Bayesian inverse problem is formulated as finding $\theta \in \R^2$
from $y$ given by
\begin{equation}
    y = \G(\theta) + \eta,  \qquad 
    \theta = \begin{bmatrix}
        \theta_{(1)} \\ 
        \theta_{(2)}
    \end{bmatrix}
    \quad \textrm{and} \quad \eta\sim\N(0, 0.1^2\I).
\end{equation}
The observations comprise pointwise measurements of $p$ and we consider
well-determined and under-determined cases:
\begin{itemize}
\item Well-determined system
    \begin{align}
        \G(\theta) = \begin{bmatrix}
            p(0.25,\theta)\\
            p(0.75,\theta)
        \end{bmatrix}  
        \qquad y = \begin{bmatrix} 
        27.5 \\ 
        79.7
        \end{bmatrix}
        \qquad
        \rho_{\rm prior} \sim \N\Bigl(\begin{bmatrix}
        0\\
        100
        \end{bmatrix}, \I\Bigr).
    \end{align}
    \item Under-determined system
    the observations 
    \begin{align}
    \G(\theta) = p(0.25,\theta) \qquad y = 27.5 \qquad \rho_{\rm prior} \sim \N\Bigl(\begin{bmatrix}
        0\\
        100
        \end{bmatrix}, \I\Bigr)
    \end{align}
\end{itemize}

The reference posterior distribution is approximated by the random walk Metropolis algorithm with a step size $1.0$ and  $5\times10^6$ samples (with a $10^6$ sample burn-in period).
We compare the  UKI-1~($J=4$), UKI-2~($J=5$), EKI, EAKI,  and ETKI applied
to (\cref{eq:evolve,eq:obser,eq:hyperparameters}), iterative Kalman filters, including IUKF-1~($J=4$), IUKF-2~($J=5$), IEnKF, IEAKF,  IETKF
all with $J = 50$ ensemble members, applied to \eqref{eq:KL_dynamics},
and the EKS and CBS methods, also  with $J = 50$.
All algorithms are initialized at the prior distribution.
The iterative Kalman filters are discretized with $\Dt = \frac{1}{30}$, and further correction~(See \ref{appendix:Gaussian_Init}) is applied on the initial ensemble members for the exactness of the initialization.

\begin{figure}[ht]
\centering
    \includegraphics[width=0.98\textwidth]{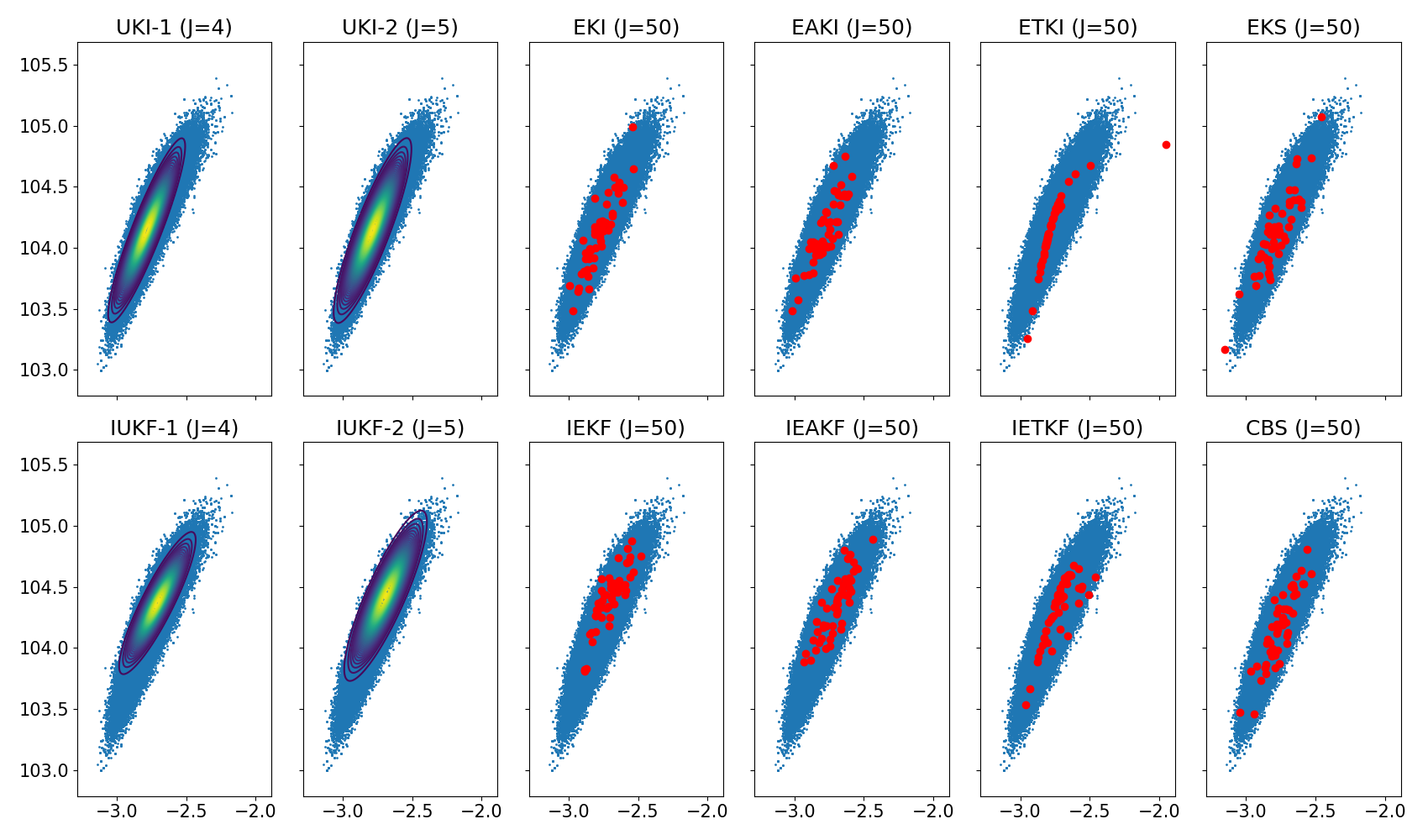}
    \caption{Nonlinear 2-parameter model problem: posterior distribution approximated at the $30^{th}$ iteration for the well-determined system. Blue dots represent the reference posterior distribution obtained by MCMC. x-axis is for $\theta_{(1)}$ and y-axis is for $\theta_{(2)}$.}
    \label{fig:nonlinear-distribution-well}
\end{figure}

Posterior distribution approximations obtained by different algorithms,
all at the $30^{th}$ iteration, are depicted in~\cref{fig:nonlinear-distribution-well,fig:nonlinear-distribution-under}.
Two common qualitative themes stand out from these figures: the iterative
methods based on coupling/transport have difficulty covering the true posterior
spread, especially in the under-determined case, when compared with the new
methodologies based on our novel mean-field dynamical system; and application of ensemble transform methods in either coupling/transport or mean-field
dynamical system suffers from a form of collapse. The first point may be
seen quantitatively; the second does not show up so much quantitatively because
collapse is in a direction in which there is less posterior spread.
We now turn to quantitative comparisons. Again, because we use the same number of steps for all algorithms, and commensurate numbers of particles, the evaluation cost of all the methods studied are comparable; the size of the error discriminates between them.

\begin{figure}[ht]
\centering
    \includegraphics[width=0.8\textwidth]{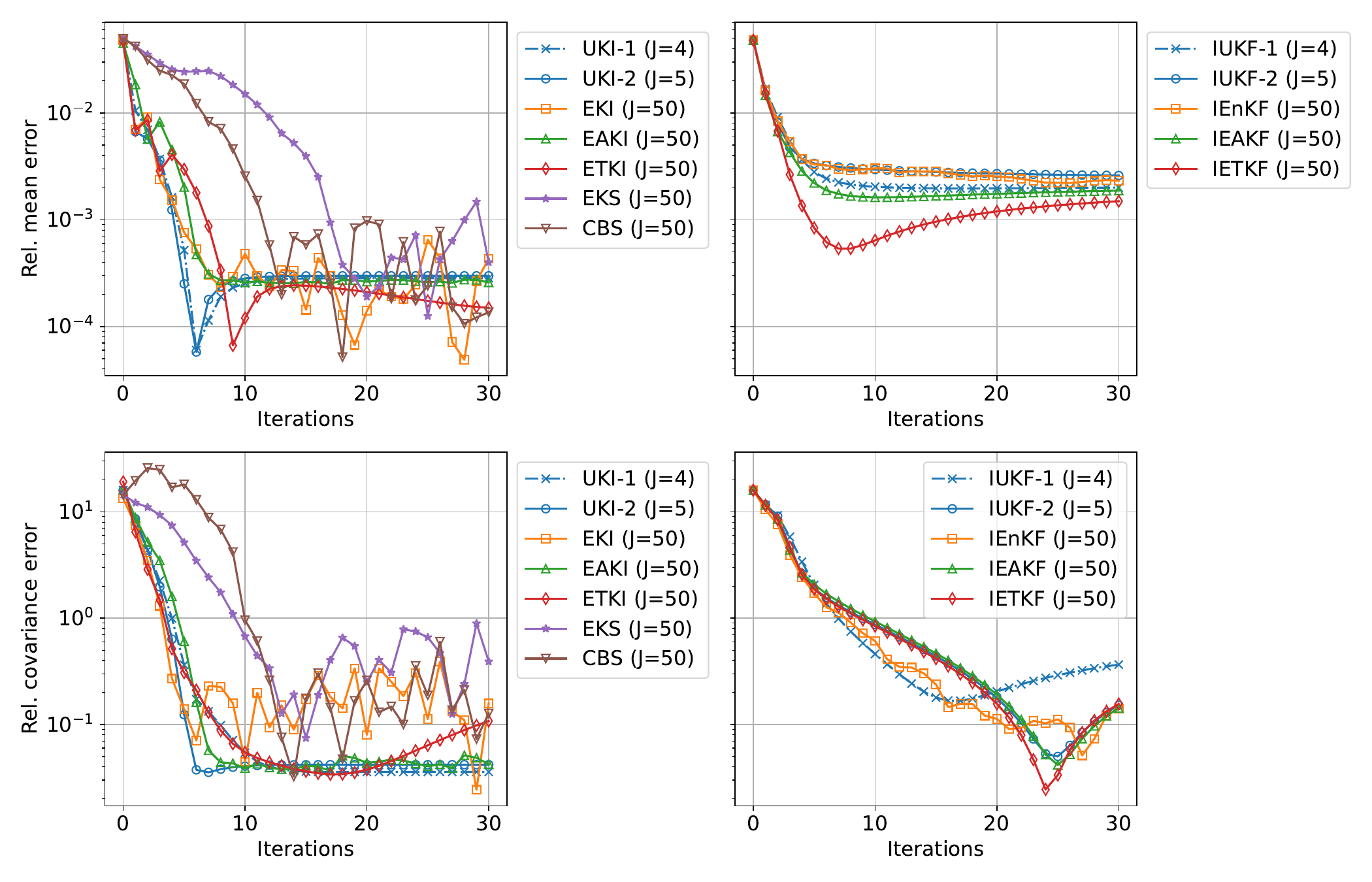}
    \caption{Nonlinear 2-parameter model problem: convergence of posterior mean (top) and posterior covariance (bottom) for the well-determined system.}
    \label{fig:nonlinear-converge-well}
\end{figure}

\begin{figure}[ht]
\centering
    \includegraphics[width=0.8\textwidth]{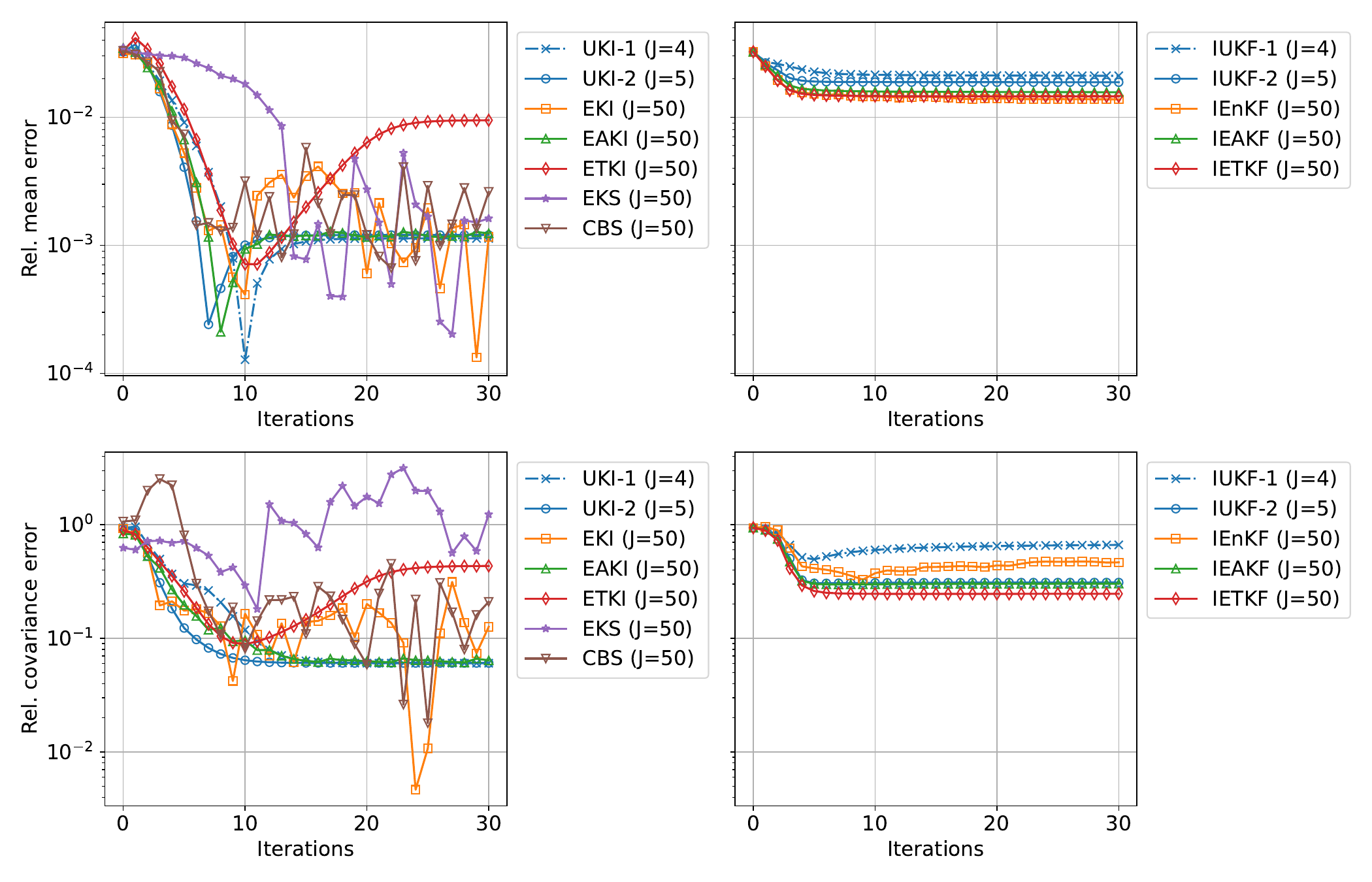}
    \caption{Nonlinear 2-parameter model problem: convergence of posterior mean (top) and posterior covariance (bottom) for the under-determined system.}
    \label{fig:nonlinear-converge-under}
\end{figure}

The convergence of posterior mean and posterior covariance are reported in~\cref{fig:nonlinear-converge-well,fig:nonlinear-converge-under}. For both scenarios, UKI-1, UKI-2, and EAKI converge exponentially fast at the beginning and then flatten out, since the posterior is not Gaussian
\cite{ernst2015analysis}.
The ETKI suffers from divergence for the under-determined scenario, and for this test, ETKI is less robust compared with UKI and EAKI.
As in the linear 2-parameter model problems~(see Subsection~\ref{ssec:app:add}), EKI, EKS, CBS, and IEnKF suffer from random noise and finite ensemble sizes.
Moreover, Kalman inversions, especially UKI and EAKI, outperform iterative Kalman filters, as measured by accuracy for commensurate cost,
for these nonlinear tests.

\begin{figure}[ht]
\centering
    \includegraphics[width=0.98\textwidth]{Figs/Elliptic-under-determined-distribution.png}
    \caption{Nonlinear 2-parameter model problem: posterior distribution approximated at the $30^{th}$ iteration for the under-determined system. Blue dots represent the reference posterior distribution obtained by MCMC. x-axis is for $\theta_{(1)}$ and y-axis is for $\theta_{(2)}$.}
    \label{fig:nonlinear-distribution-under}
\end{figure}

\subsection{Hilbert Matrix Problem}
\label{sec:app-Hilbert}

We define the Hilbert matrix $G \in \R^{N_{\theta}\times N_{\theta}}$ by its entries
\begin{align}
    G_{i,j} = \frac{1}{i+j-1},
\end{align}
with $N_{\theta} = 100.$
We consider the inverse problem
\begin{align}
y = G\mathds{1} \qquad \eta \sim \N(0, 0.1^2\I) \qquad \rho_{\rm prior} \sim \N(0, \I).
\end{align}

We no longer study iterated Kalman methods arising from coupling/transport
as the preceding examples show that they are inefficient. Furthermore
EKI, EKS, and CBS
do not converge and  suffer from random noise and/or finite ensemble sizes;
these results are not shown.
Instead, we focus on the effect of the ensemble size on EAKI and ETKI,
comparing with UKI. To be concrete, we apply EAKI and ETKI with $J=N_{\theta}$, $N_{\theta}+1$ and $500$, and UKI-1 and UKI-2. Again, we initialize all algorithms at the prior distribution.

\begin{figure}[ht]
\centering
    \includegraphics[width=0.8\textwidth]{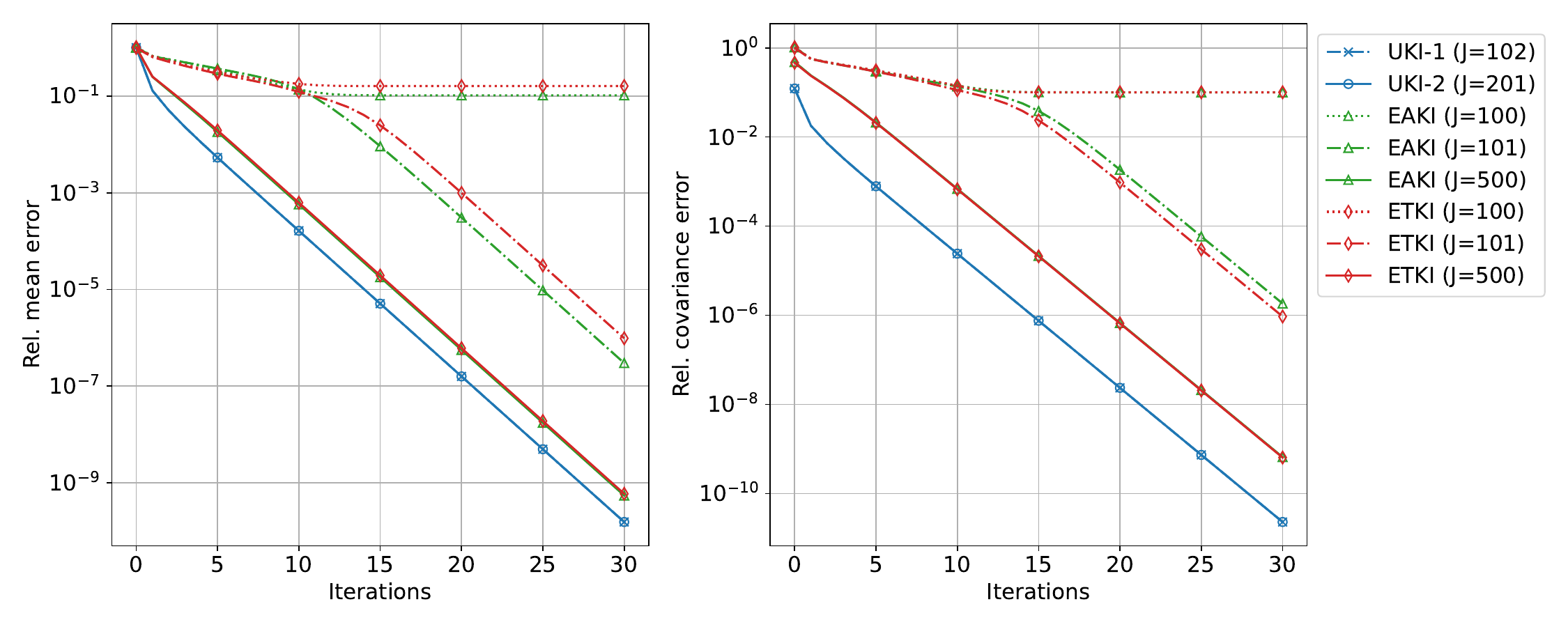}
    \caption{Hilbert matrix problem: convergence of posterior mean (left) and posterior covariance (right).}
    \label{fig:Hilbert}
\end{figure}

We compute the reference distribution analytically. The convergence of posterior mean and posterior covariance are reported in~\cref{fig:Hilbert}. 
UKI-1, UKI-2, and EAKI and ETKI with more than $N_{\theta}$ ensemble particles, converge exponentially fast. 
The relatively poor performance of EAKI and ETKI  with a smaller 
number of ensemble particles is related to \cref{th:lin_converge}, since EAKI and ETKI require at least $N_{\theta}+1$ particles to ensure that the initial covariance matrix $\Cov_0$ is strictly positive definite.

\subsection{Darcy Flow Problem}

The two-dimensional Darcy flow equation describes the pressure field $p(x)$ in a porous medium defined by a parameterized, positive permeability field $a(x,\theta)$:
\begin{equation}
    \begin{aligned}
       -\nabla \cdot (a(x, \theta) \nabla p(x)) &= f(x), \quad &&x\in D,\\
    p(x) &= 0, \quad &&x\in \partial D.         
    \end{aligned}
\end{equation}
Here the computational domain is $D = [0,1]^2$,  Dirichlet boundary conditions are applied on  $\partial D$, and $f$ defines the source of the fluid:
\begin{align}
    f(x_1, x_2) = \begin{cases}
               1000 & 0 \leq x_2 \leq \frac{4}{6}\\
               2000 & \frac{4}{6} < x_2 \leq \frac{5}{6}\\
               3000 & \frac{5}{6} < x_2 \leq 1\\
            \end{cases}. 
\end{align} 

\begin{figure}[ht]
\centering
    \includegraphics[width=0.48\textwidth]{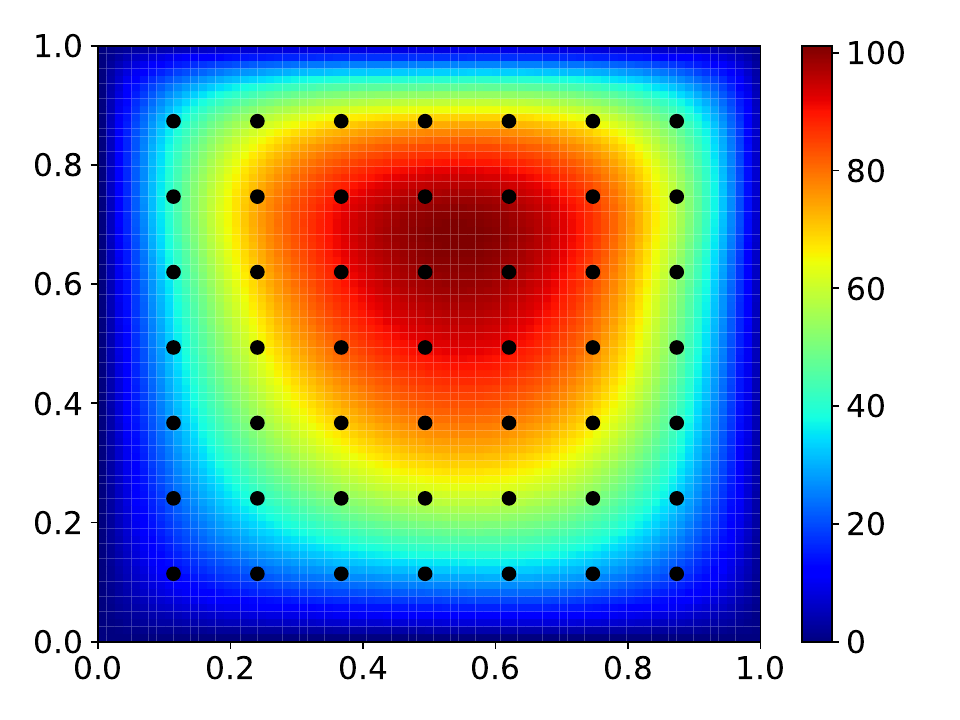}
    \caption{The reference pressure field with $49$ equidistant pointwise measurements of the Darcy flow problem.}
    \label{fig:darcy-1d-ref}
\end{figure}

The inverse problem of interest is to determine parameter $\theta$ of
the field $a(\cdot;\theta)$ from observation $y_{ref}$, which consists of pointwise measurements of the pressure value $p(\cdot)$ at $49$ equidistant points in the domain~(see~\cref{fig:darcy-1d-ref}), corrupted with observation error $\eta \sim \N(0, \I)$. We now describe how $a$ depends on $\theta$, and specify a standard Gaussian prior on $\theta.$

We write
\begin{equation}
\label{eq:KL-2d}
    \log a(x,\theta) = \sum_{l\in K} \theta_{(l)}\sqrt{\lambda_l} \psi_l(x),
\end{equation}
where $K = \Z^{+}\times\Z^{+} \setminus \{0,0\}$, and 
\begin{equation}
    \psi_l(x) = \begin{cases}
                 \sqrt{2}\cos(\pi l_1 x_1)              & l_2 = 0\\
                 \sqrt{2}\cos(\pi l_2 x_2)              & l_1 = 0\\
                 2\cos(\pi l_1 x_1)\cos(\pi l_2 x_2)    & \textrm{otherwise}\\
                 \end{cases},
                 \qquad \lambda_l = (\pi^2 |l|^2 + \tau^2)^{-d}
\end{equation}
and $\theta_{(l)} \sim \N(0,1)$ i.i.d. The expansion~\cref{eq:KL-2d} can be rewritten as a sum over $\Z^{+}$ rather than a lattice: 
\begin{equation}
\label{eq:KL-1d}
    \log a(x,\theta) = \sum_{k\in \Z^{+}} \theta_{(k)}\sqrt{\lambda_k} \psi_k(x),
\end{equation}
where the eigenvalues $\lambda_k$ are in descending order.
We note that these considerations amount to assuming
that $\log a(x, \theta)$ is a mean zero Gaussian random field
with covariance
\begin{equation}
    \mathsf{C} = (-\Delta + \tau^2 )^{-d},
\end{equation}
with $-\Delta$ the Laplacian on $D$ subject to homogeneous Neumann boundary conditions on the space of spatial-mean zero functions; hyperparameter 
$\tau$ denotes the inverse length scale of the random field and hyperparameter $d$ determines its Sobolev and H\"older regularity, which is $d-1$
in our two dimensional setting \cite{dashti2013bayesian}.

In this work, we take $\tau=3$ and $d=2.$
In practice, we truncate the sum \eqref{eq:KL-1d} to $N_\theta$ terms, 
based on the largest $N_\theta$ eigenvalues, and hence $\theta\in\R^{N_\theta}$. The forward problem is solved by a finite difference method on a $80 \times 80$ grid. 
To create the data $y_{ref}$ referred to above, we generate a truth random field $\log a(x,\theta_{ref})$ with $N_{\theta} = 128$ and $\theta_{ref} \sim \N(0, \I^{128})$.

The benchmark posterior distribution is approximated by the preconditioned Crank–Nicolson algorithm with $2\times 10^6$ samples~(with a $5\times10^5$ sample burn-in period) with the step size $0.04$. Since the preceding
examples have shown the benefits of using UKI, EAKI and ETKI over all other
methods considered, we compare only these approaches with the benchmark.
Specifically, we apply UKI-1, UKI-2, and  EAKI and ETKI with $J=N_{\theta}+2$, again initialized at the prior distribution.

\begin{figure}[ht]
\centering
    \includegraphics[width=0.9\textwidth]{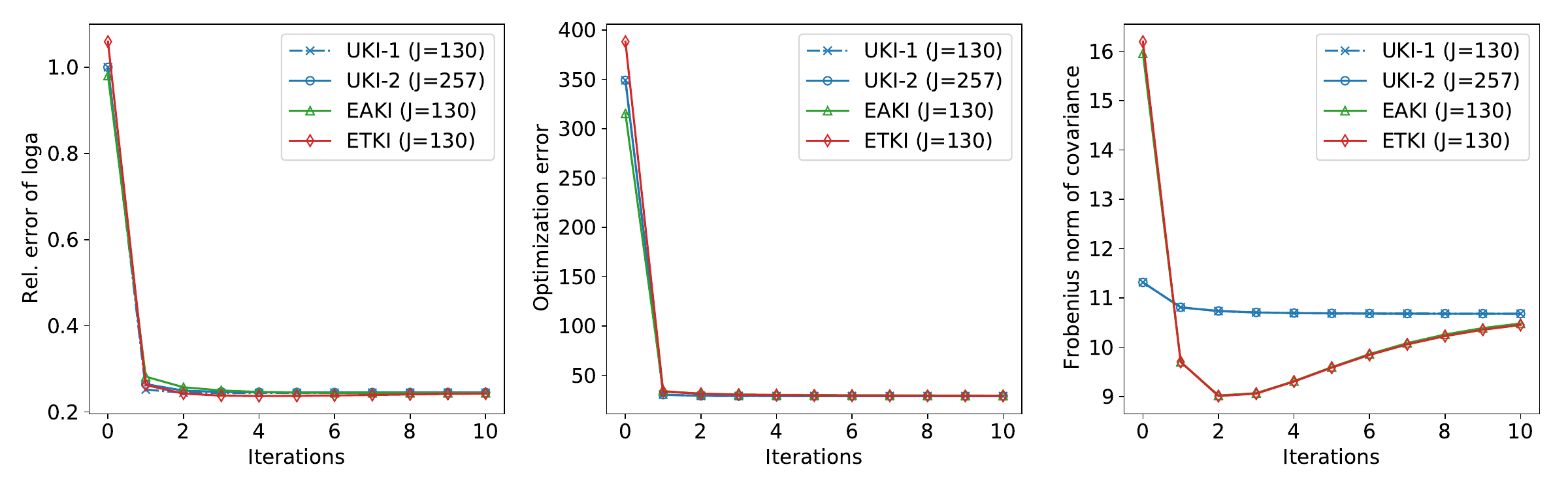}
    \caption{
    Darcy flow problem: the relative error of the permeability field, the optimization error $\displaystyle \frac{1}{2}\lVert\Sigma_{\eta}^{-\frac{1}{2}} (y_{obs} - \py_n)\rVert^2 + \frac{1}{2}\lVert\Sigma_{0}^{-\frac{1}{2}} (\mean_{n} - r_0)\rVert^2$ and the Frobenius norm $\lVert\Cov_n\rVert$ (from left to right).}
    \label{fig:Darcy-converge}
\end{figure}

The convergence of the relative $L_2$ error of the mean of the $\log a$ field, the optimization errors, and the Frobenius norm of the estimated posterior covariance, as the iteration progresses, are depicted in~\cref{fig:Darcy-converge}. This clearly shows that all four Kalman inversion techniques converge within $10$ iterations.

\begin{figure}[ht]
\centering
    \includegraphics[width=0.98\textwidth]{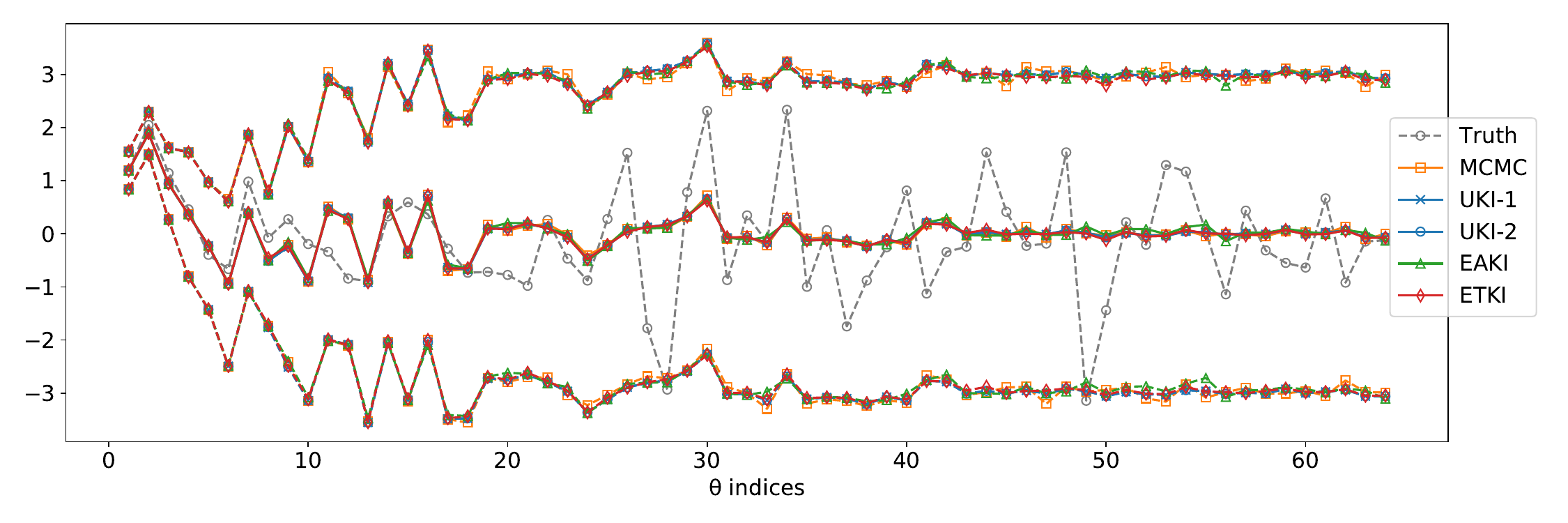}
    \caption{The estimated KL expansion parameters $\theta_{(i)}$ and the associated 3-$\sigma$ confidence intervals obtained by UKI-1~($J=130$), UKI-2~($J=257$), EAKI~($J=130$), ETKI~($J=130$) and MCMC for the Darcy flow problem.}
    \label{fig:darcy-2d-uq}
\end{figure}

\begin{figure}[ht]
\centering
    \includegraphics[width=0.98\textwidth]{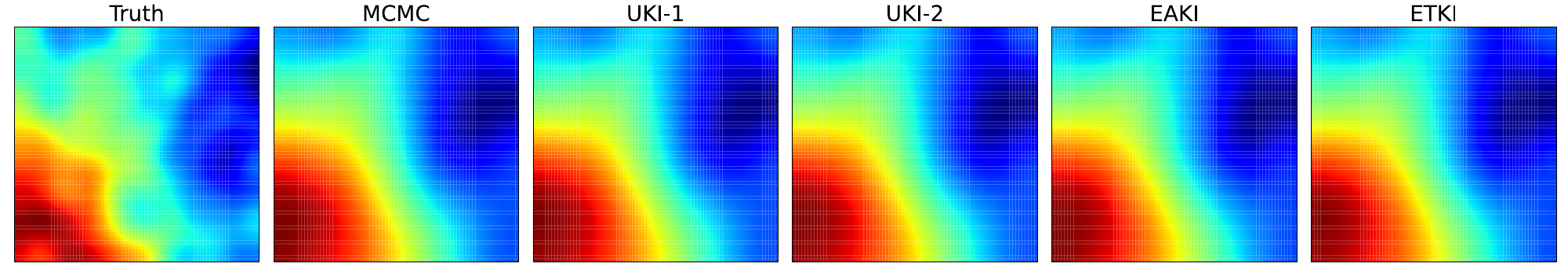}
    \caption{The truth log permeability field $\log a$, and log permeability fields obtained by MCMC, UKI-1~($J=130$), UKI-2~($J=257$), EAKI~($J=130$), ETKI~($J=130$) (left to right).}
    \label{fig:darcy-logk}
\end{figure}

\cref{fig:darcy-2d-uq} shows the properties of the converged
posterior distribution, after the $10^{th}$ iteration, comparing them
with MCMC and with the truth (referred to as ``Truth''). The information
is broken down according to recovery of the $\{\theta_{(i)}\}$,
visualizing only the first $64$ modes, since the statistical estimates
of other modes obtained by MCMC and by our Kalman inversion methodologies
are close to the prior $\N(0,1)$ -- the data does not inform them.  
We first note that the truth values lie in the confidence intervals
determined by MCMC, with high probabilities. Secondly, we note that
all four Kalman methods reproduce the posterior mean and confidence intervals
computed by MCMC accurately. The estimated log-permeability fields $\log a$ and the truth  are depicted in \cref{fig:darcy-logk}. 
The mean estimations obtained by the MCMC and these Kalman inversions match well, and they both capture the main feature of the truth log-permeability field.

\begin{figure}[ht]
\centering
    \includegraphics[width=0.98\textwidth]{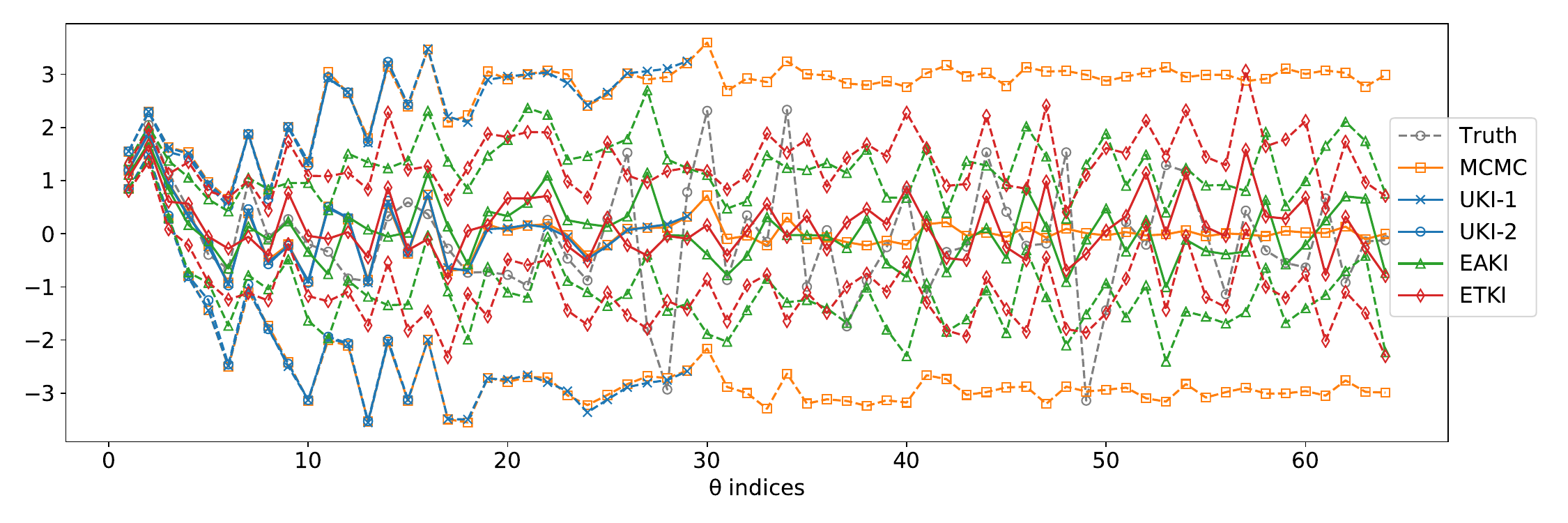}
    \caption{The estimated KL expansion parameters $\theta_{(i)}$ and the associated 3-$\sigma$ confidence intervals obtained by UKI-1, UKI-2, EAKI, ETKI with $J=31$ ensemble members and MCMC for the Darcy flow problem.}
    \label{fig:darcy-2d-uq-LR}
\end{figure}

\begin{remark}
In practice, for many inversion problems for fields,
the realistic number of ensemble members is much smaller than the dimension of the state space. To probe this setting, we repeat the test by using UKI-1, UKI-2, ETKI, and EAKI with $J=31$ ensemble members; for UKI-1 and UKI-2, we invert
for the first $29$ and $15$ coefficients of $\{\theta_{(k)}\}$,
respectively, and for ETKI and EAKI we invert for all $128$ coefficients. 
The estimated KL expansion parameters $\{\theta_{(i)}\}$ for the log-permeability field and the associated  3-$\sigma$ confidence intervals obtained by MCMC, and different Kalman inversions at the $10^{th}$ iteration, are depicted in~\cref{fig:darcy-2d-uq-LR}. The mean and standard deviation of the coefficients associated with these dominant modes obtained by both UKIs match well with those obtained by MCMC. The results indicate that
the ``truncate then invert'' strategy used by UKIs outperforms the ``direct inversion'' strategy, used here by EAKI and ETKI, when only small ensemble
numbers are feasible.
\end{remark}

\subsection{Idealized Global Climate Model}

\label{sec:app:GCM}
Finally, we consider using low-fidelity model techniques to speed up an idealized global climate model inverse problem. The model is based on the 3D Navier-Stokes equations, making the hydrostatic and shallow-atmosphere approximations common in atmospheric modeling. 
Specifically, we test on the notable Held-Suarez test case~\cite{held1994proposal}, in which a detailed radiative transfer model is replaced by Newtonian relaxation of temperatures toward a prescribed ``radiative equilibrium'' $T_{\mathrm{eq}}(\phi, p)$ that varies with latitude $\phi$ and pressure $p$. \dzh{Specifically, the thermodynamic equation for temperature $T$ 
\begin{align}
\frac{D T}{\partial t} - \frac{R T \omega}{C_p p} = Q
\end{align}
(including advective and pressure work terms)} contains a diabatic heat source 
\begin{align}
Q = -k_T(\phi, p, p_s) \bigl(T - T_{\mathrm{eq}}(\phi, p)\bigr),
\end{align}
with relaxation coefficient (inverse relaxation time)
\begin{align}
k_T = k_a + (k_{s} - k_a)\max\Bigl(0, \frac{\sigma - \sigma_b}{1 - \sigma_b}\Bigr)\cos^4\phi.
\end{align}
Here, $\sigma = p/p_s$, pressure $p$ normalized by surface pressure $p_s$, is the vertical coordinate of the model, and 
\begin{align}
T_{\mathrm{eq}} = \max\Bigl\{200K, \Bigl[315K - \Delta T_y \sin^2\phi - \Delta\theta_z\log\Bigl(\frac{p}{p_0}\Bigr)\cos^2\phi\Bigr]\Bigl(\frac{p}{p_0}\Bigr)^{\kappa}\Bigr\}
\end{align}
is the equilibrium temperature profile ($p_0 = 10^5~\mathrm{Pa}$ is a reference surface pressure and $\kappa = 2/7$ is the adiabatic exponent).

The inverse problem of interest here is to
determine the parameters $(k_a,\ k_{s},\ \Delta T_y,\ \Delta\theta_z)$ 
from statistical averages of the temperature field $T$.
We impose the following constraints: 
$$0\  \mathrm{day}^{-1}<k_a < 1\  \mathrm{day}^{-1}, \quad   0\  \mathrm{day}^{-1}<k_s <1\  \mathrm{day}^{-1}, \quad   0\  \mathrm{K} < \Delta T_y < 100\  \mathrm{K}, \quad   0\  \mathrm{K}<\Delta\theta_z< 50\  \mathrm{K}.$$
The inverse problem is formed as follows~\cite{UKI},
\begin{equation}
    y = \G(\theta) + \eta \quad \mathrm{with} \quad \G(\theta)=\overline{T}(\phi, \sigma)
\end{equation}
with the parameter transformation 
\begin{equation}
    \theta: (k_a, k_s, \Delta T_y, \Delta\theta_z) = \Big(
    \frac{1}{1+\exp(\theta^{(1)})},\ 
    \frac{1}{1+\exp(\theta^{(2)})},\
    \frac{100}{1+\exp(\theta^{(3)})},\
    \frac{50}{1+\exp(\theta^{(4)})}
    \Big)
\end{equation}
enforcing the constraints.
The observation mapping $\G$ is defined by mapping from the unknown $\theta$ 
to the $200$-day zonal mean of the temperature~($\overline{T}$) as a function of latitude ($\phi$) and height ($\sigma$), after an initial spin-up 
of $200$ days.

\begin{figure}[ht]
\centering
\includegraphics[width=0.48\textwidth]{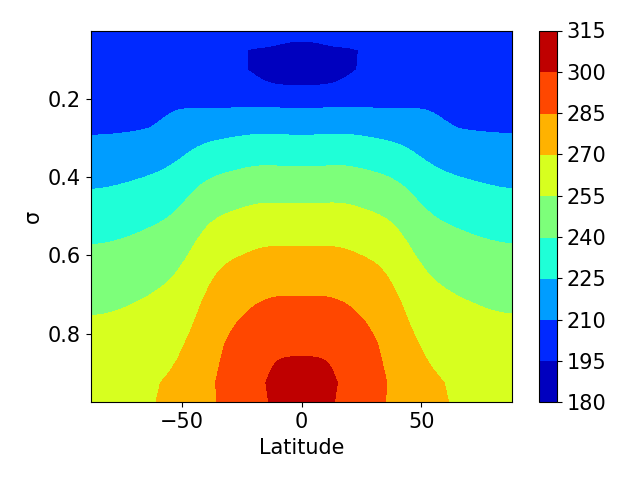}
\includegraphics[width=0.48\textwidth]{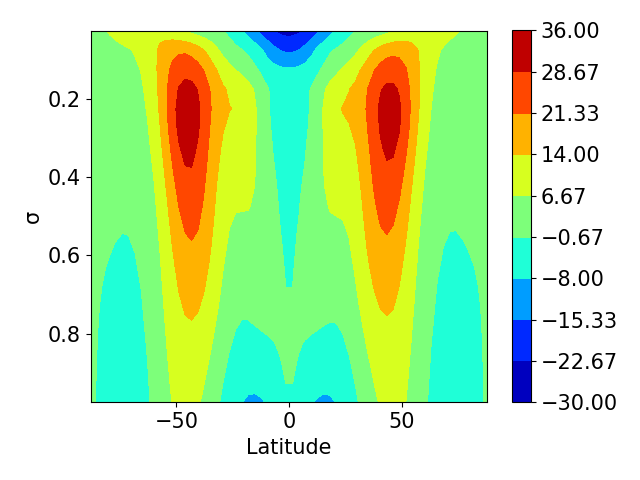}\\
\includegraphics[width=0.48\textwidth]{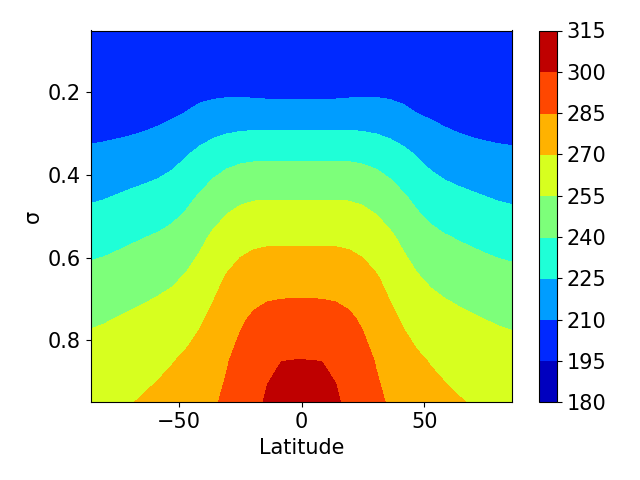}
\includegraphics[width=0.48\textwidth]{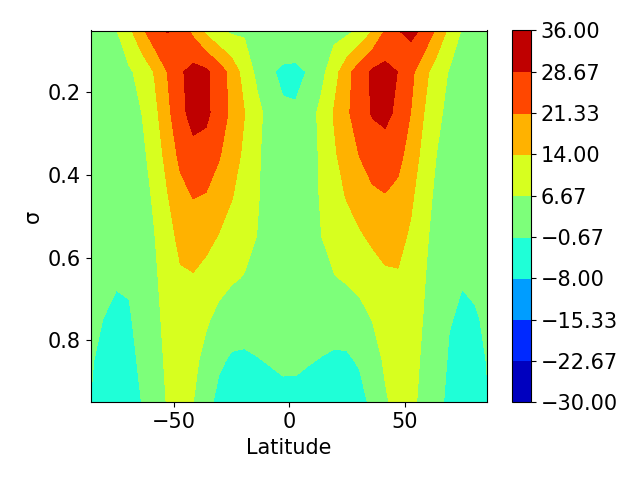}
\caption{Zonal mean temperature~(left) and zonal wind velocity~(right) obtained with the T42 grid~(top) and the T21 grid~(bottom).}
\label{fig:GCM-sol}
\end{figure}

Default parameters used to generate the data in our simulation study are
\[
k_a = (40\ \mathrm{day})^{-1},\qquad  
k_{s} = (4\ \mathrm{day})^{-1}, \qquad
\Delta T_y = 60\ \mathrm{K}, \qquad
\Delta\theta_z = 10\ \mathrm{K}.
\]
For the numerical simulations, we use the spectral transform method in the horizontal, with T42 spectral resolution~(triangular truncation at wavenumber 42, with $64 \times 128$ points on the latitude-longitude transform grid); we use 20 vertical levels equally spaced in $\sigma$. With the default parameters, the model produces an Earth-like zonal-mean circulation, albeit without moisture or precipitation.  
The truth observation is the 1000-day zonal mean of the temperature~(see \cref{fig:GCM-sol}-top-left), after an initial spin-up, also of 200 days, to eliminate the influence of the initial condition. Because the truth observations come from an average 5 times as long as the observation window used for parameter learning, the chaotic internal variability of the model introduces noise in the observations.

To perform the inversion, we set the prior $\displaystyle \rho_{\rm prior} \sim \N(0,\ 10^2\I)$. Within the algorithm, we assume that the observation error satisfies $\eta \sim \N(0\ \textrm{K}, 3^2\I\ \textrm{K}^2)$. 
All these Kalman inversions are initialized with $\theta_0 \sim \N(0,\ 0.1^2\I)$, since initializing at the prior leads to unstable simulations at the first iteration. 
The bi-fidelity approach discussed in Subsection~\ref{subs:bi-fidelity} is applied to speed up both UKI-1 and UKI-2. These $J-1$ forward model evaluations are computed on a T21 grid~(triangular truncation at wavenumber 21, with $32 \times 64$ points on the latitude-longitude transform grid) with 10 vertical levels equally spaced in $\sigma$~(twice coarser in all three directions). They are abbreviated as UKI-1-BF and UKI-2-BF. 
The computational cost of the high-fidelity~(T42) and low-fidelity~(T21) models are about $4$-CPU hour and $0.5$-CPU hour, and therefore the bi-fidelity approach effectively reduces CPU costs.
The 1000-day zonal mean of the temperature and velocity predicted by the low-resolution model with the truth parameters are shown in~\cref{fig:GCM-sol}-bottom. It is worth mentioning there are significant discrepancies comparing with results computed on the T42 grid~(\cref{fig:GCM-sol}-top). Whether these would be tolerable will depend on the
use to which the posterior inference is put.

\begin{figure}[ht]
\centering
\includegraphics[width=0.45\textwidth]{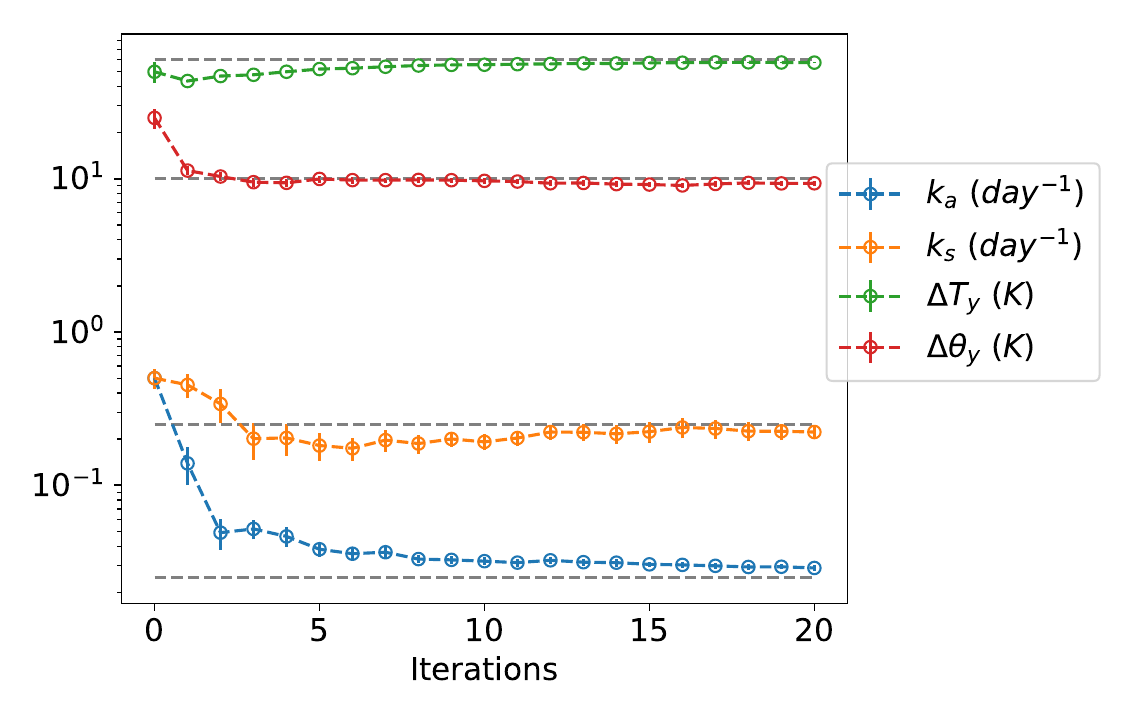}
\includegraphics[width=0.45\textwidth]{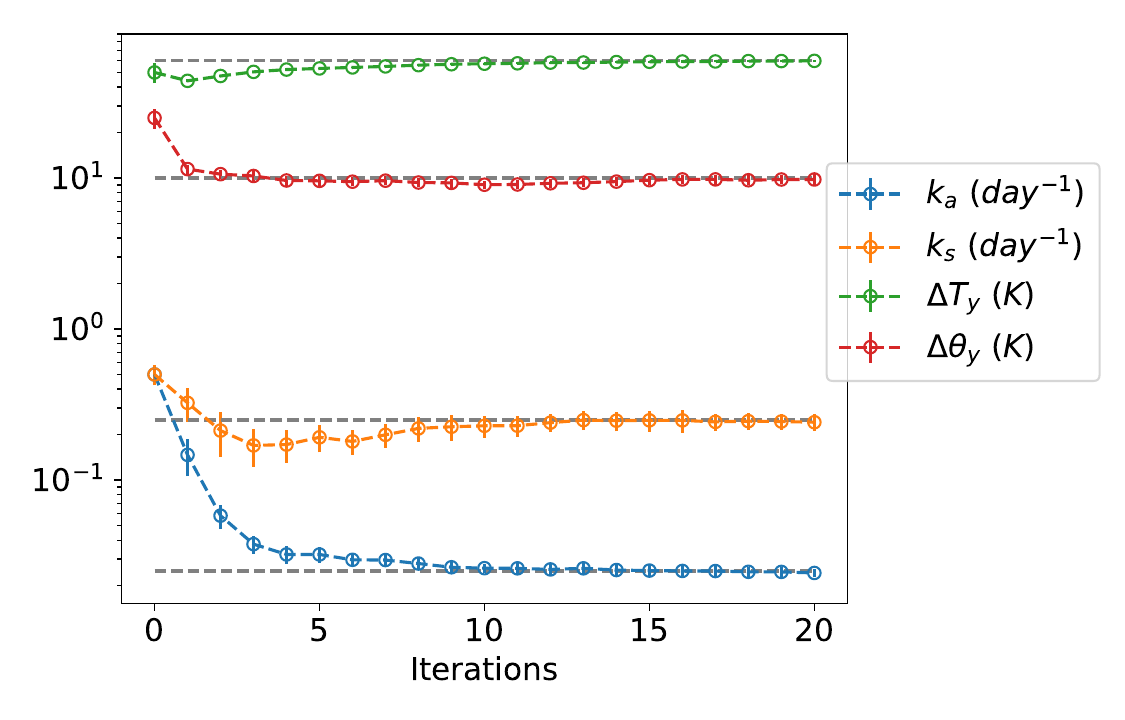}
\includegraphics[width=0.45\textwidth]{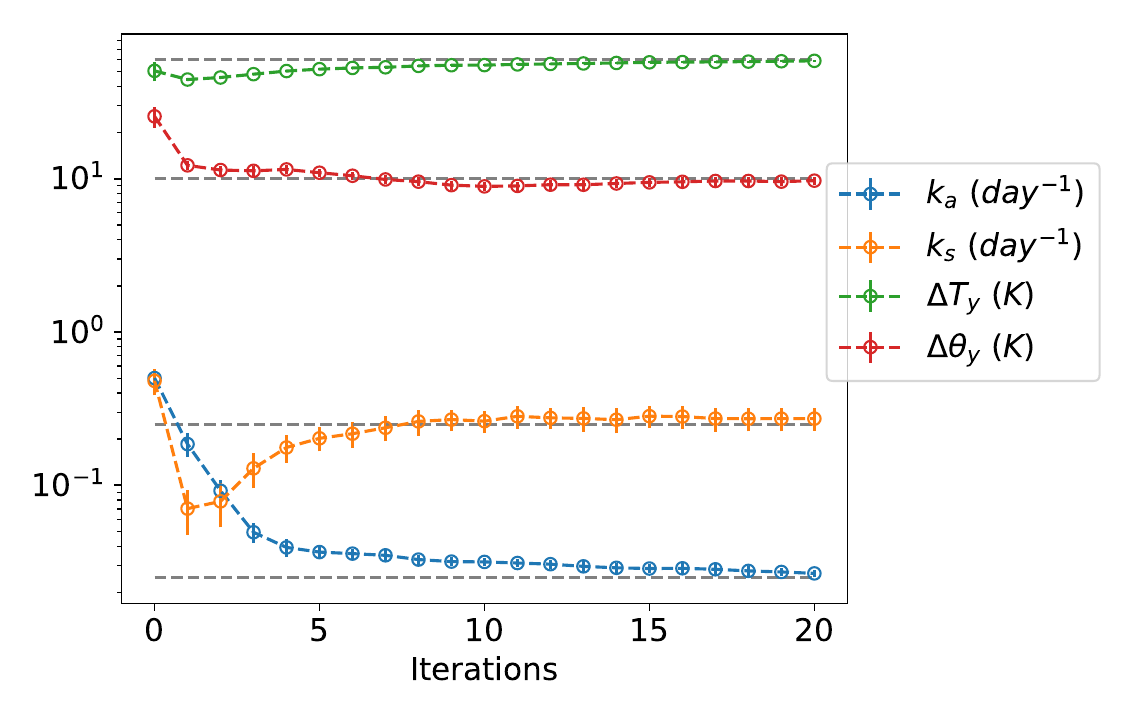}
\includegraphics[width=0.45\textwidth]{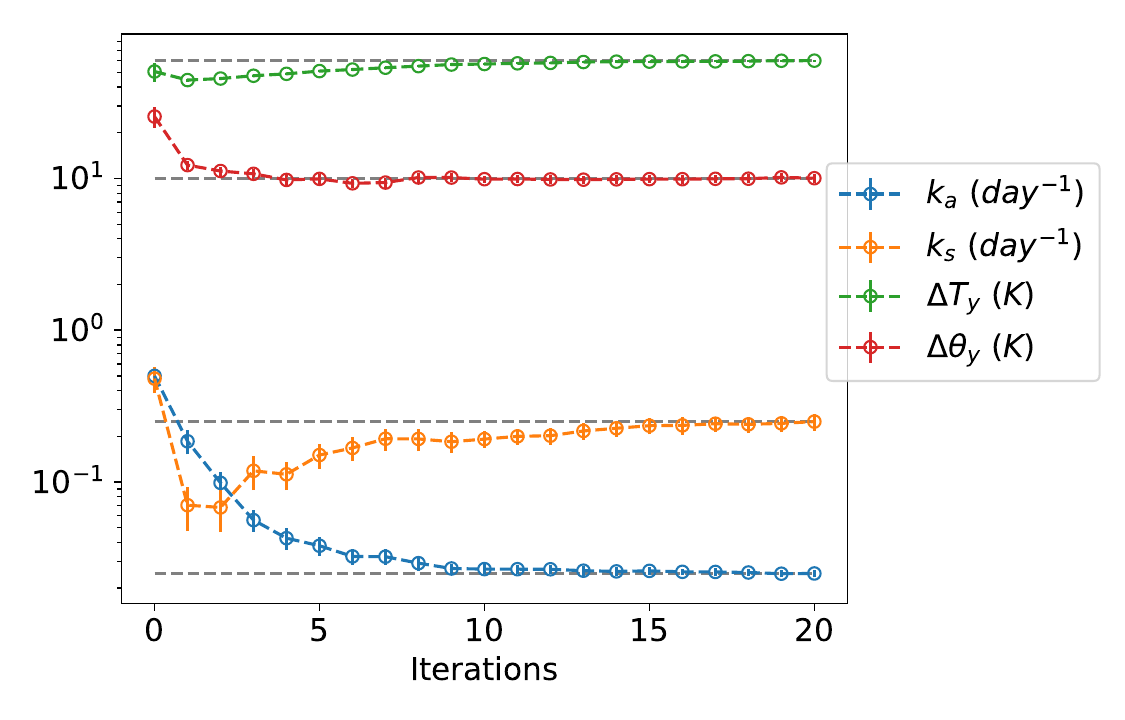}
\includegraphics[width=0.45\textwidth]{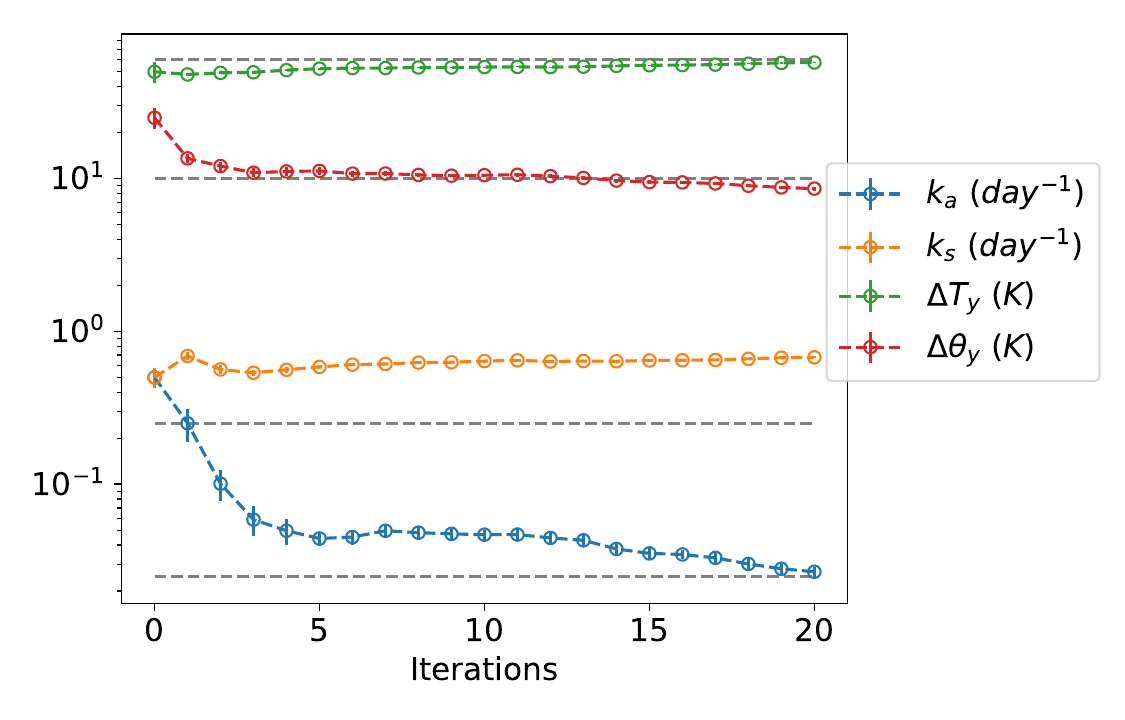}
\includegraphics[width=0.45\textwidth]{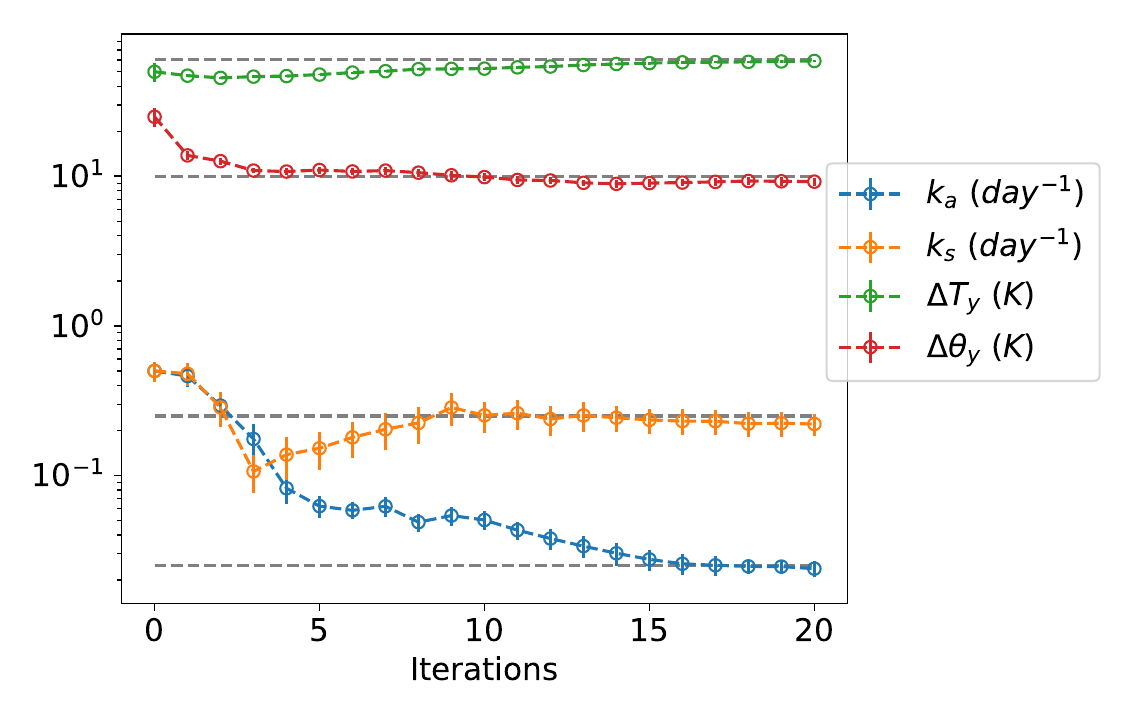}
\caption{Convergence of the idealized global climate model inverse problem. 
True parameter values are represented by dashed grey lines.
Top: UKI-1~($J=6$) and UKI-2~($J=9$); 
Middle: EAKI~(J=$9$) and ETKI~($J=9$);
Bottom: UKI-1-BF~($J=6$) and UKI-2-BF~($J=9$).
}
\label{fig:GCM-obj}
\end{figure}

The estimated parameters and associated $3-\sigma$ confidence intervals for each component at each iteration are depicted in \cref{fig:GCM-obj}.
Since the prior covariance is large and the problem is over-determined it is natural to expect that
the posterior mean should be close to the true parameters. Indeed all the
different Kalman inversions, except UKI-1-BF~($J=6$), do indeed converge to the true parameters.

\section{Conclusion}
\label{sec:C}
Kalman-based inversion has been widely used to construct derivative-free optimization and sampling methods for nonlinear inverse problems.
In this paper, we developed new Kalman-based inversion methods, for Bayesian inference and uncertainty quantification, which build on the work in both
optimization and sampling.
We  propose  a  new  method  for  Bayesian  inference  based  on filtering  a novel  mean-field dynamical system subject to partial noisy observations,
and which depends on the law of its own filtering distribution, 
together with application of the Kalman methodology. 
Theoretical guarantees are presented: for linear inverse problems, the mean and covariance obtained by the method converge exponentially fast to the posterior mean and covariance.
For nonlinear inverse problems, numerical studies indicate the method delivers an excellent approximation of the posterior distribution for problems which
are not too far from Gaussian.
The methods are shown to be superior to existing coupling/transport
methods, under the umbrella of iterative Kalman methods; and deterministic
rather than stochastic implementations of Kalman methodology are found to
be favorable.
We further propose several simple strategies, including low-rank approximation and a bi-fidelity approach, to reduce the computational and memory cost of the
proposed methodology.

There are numerous directions for future study in this area. On the theoretical
side it would be of value to obtain theoretical guarantees concerning the
accuracy of the methodology, when applied to nonlinear inverse problems, which are
close to Gaussian. On the methodological side, it would be of particular interest
to systematically quantify model-form error by using this framework, rather than to
assume perfect the model scenario as we have done here.

\paragraph{Acknowledgments} AMS and DZH are supported by the generosity of Eric and Wendy Schmidt by recommendation of the Schmidt Futures program; AMS is also
supported by the Office of Naval Research (ONR) through grant N00014-17-1-2079. 
JH is supported by the Simons Foundation as a Junior Fellow at New York University. SR is supported by Deutsche Forschungsgemeinschaft (DFG) - Project-ID 318763901 - SFB1294.

\appendix
\section{Continuous Time Limit}
\label{appendix:CTL}

To derive a continuous-time limit
of the novel mean-field dynamical system (\cref{eq:evolve,eq:obser,eq:hyperparameters}) we define
$\tau_n=n\Dt$ and define $\{z_n\}$
by $x_n=\Dt^{-1}(z_{n+1}-z_n).$ Then  $\theta_n \approx \theta(\tau_n)$
and $z_n \approx z(\tau_n)$.
Let $W \in \R^{N_{\theta}}$ and $B\in \R^{N_y + N_{\theta}}$ be standard unit Brownian motions. Letting $\Dt \rightarrow 0$ in 
(\cref{eq:evolve,eq:obser,eq:hyperparameters}), we then obtain
\begin{equation}
\begin{aligned}
    &\dot{\theta} =  C^{\frac{1}{2}} \dot{W},\\
    &\dot{z} = \bG(\theta) + \begin{bmatrix}
       \Sigma_{\eta} & 0\\
        0 & \Sigma_{0} 
        \end{bmatrix}^{\frac{1}{2}}\dot{B}.
\end{aligned}
\end{equation}
We are interested in the filtering problem of finding the
distribution of $\theta(\tau)$ given $\{z(s)\}_{s=0}^{\tau}$ and then
evaluating this distribution in the setting \dzh{where
$\dot{z}(s) \equiv x$ defined in \eqref{eq:data}.}

Under a similar limiting process,
the Gaussian approximation algorithm defined by \eqref{eq:sds} becomes
\begin{equation}
\label{eq:ct-limit2}
\begin{aligned}
\dot{m} &=  \pCov^{\theta \by} \begin{bmatrix}
        \Sigma_{\eta} & 0\\
        0 & \Sigma_{0} 
        \end{bmatrix}^{-1} (x - \bpy),\\
\dot{C}&= C -  \pCov^{\theta \by}\begin{bmatrix}
        \Sigma_{\eta} & 0\\
        0 & \Sigma_{0} 
        \end{bmatrix}^{-1} {\pCov^{\theta \by}}{}^{T}, 
\end{aligned}
\end{equation}
where 
$\bpy = \E(\bG(\theta))$, 
$\pCov^{\theta \by}=\E\Bigl(\bigl(\theta-m\bigr)\otimes\bigl(\bG(\theta)-\E \bG(\theta)\bigr)\Bigr)$,
and expectation $\E$ is with respect to the distribution $\theta \sim \N(m,C)$. 
Define $\hat\Phi_R(m, C) = (x - \bpy)^T\begin{bmatrix}
        \Sigma_{\eta} & 0\\
        0 & \Sigma_{0} 
        \end{bmatrix}^{-1} (x - \bpy)$, and note that 
        $\frac{\partial \bpy}{\partial m} = \pCov^{\theta  \by^T}\Cov^{-1}$~\cite[Lemma 2]{UKI}. It follows that \cref{eq:ct-limit2} can be rewritten as 
\begin{equation}
\label{eq:ct-limit3}
\begin{aligned}
\dot{m} &= -\frac{1}{2} \Cov \frac{\partial \hat\Phi_R}{\partial m},\\
\dot{C}^{-1}&= - C^{-1} +  \Big(\frac{\partial \bpy}{\partial m} \Big)^T\begin{bmatrix}
        \Sigma_{\eta} & 0\\
        0 & \Sigma_{0} 
        \end{bmatrix}^{-1} \Big(\frac{\partial \bpy}{\partial m} \Big).
\end{aligned}
\end{equation}
The stationary points of \cref{eq:ct-limit3} satisfy
\begin{equation}
\begin{aligned}
\frac{\partial \hat\Phi_R}{\partial m}=0 \quad \textrm{and} \quad
 C^{-1} = \Big(\frac{\partial \bpy}{\partial m} \Big)^T\begin{bmatrix}
        \Sigma_{\eta} & 0\\
        0 & \Sigma_{0} 
        \end{bmatrix}^{-1} \Big(\frac{\partial \bpy}{\partial m} \Big) \succeq \Sigma_{0}^{-1}.
\end{aligned}
\end{equation}

\begin{figure}[ht]
\centering
    \includegraphics[width=0.8\textwidth]{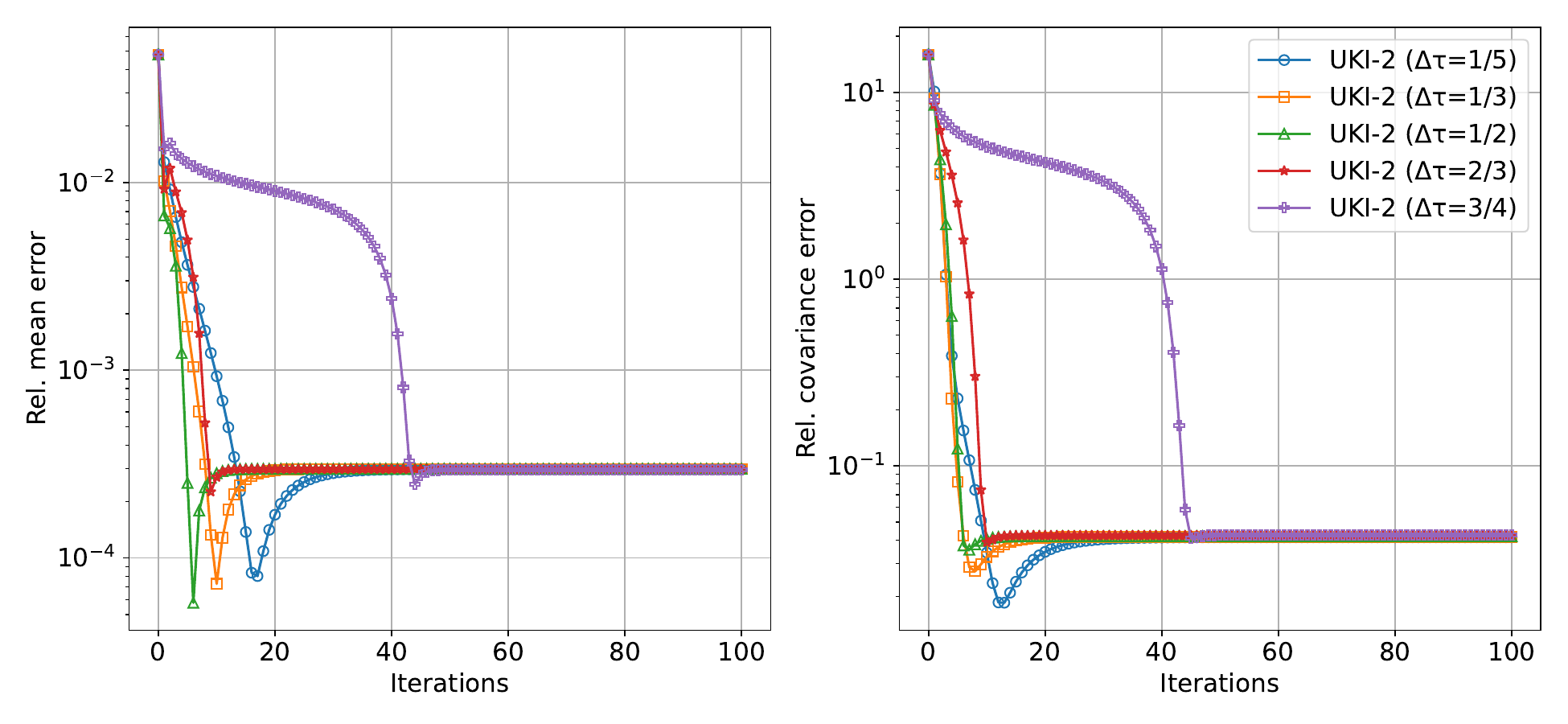}
    \includegraphics[width=0.8\textwidth]{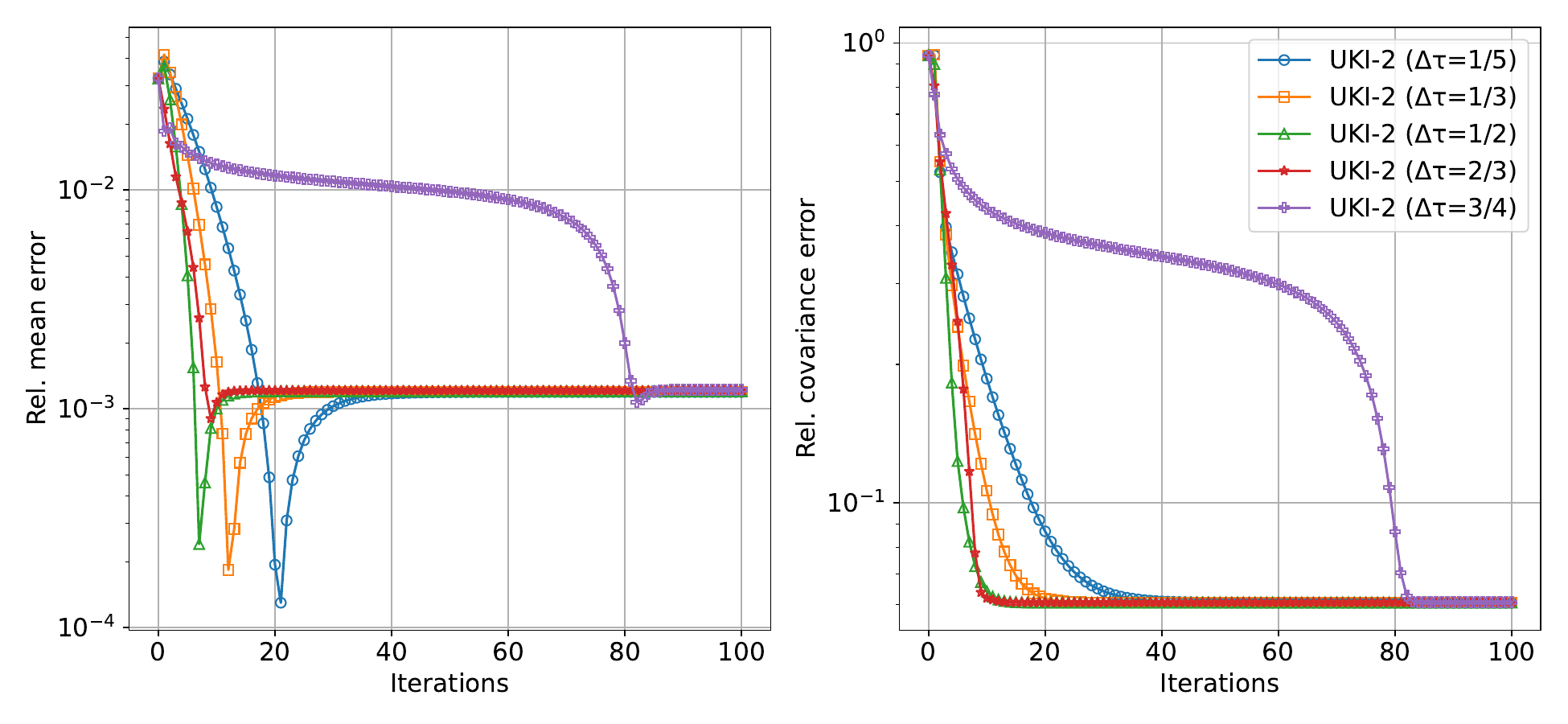}
    \caption{Nonlinear 2-parameter model problem: convergence of posterior mean (left) and posterior covariance (right) for the well-determined system (top) and the under-determined system (bottom) with different $\Dt$.}
    \label{fig:nonlinear-converge-gamma}
\end{figure}

\dzh{Although the $\Dt$ plays no role in the continuous
time limit, the original $\Dt$ parameter}
does have a marked effect on discrete-time algorithms used in practice.
To highlight this, we study the empirical effect of $\Dt$, in the context of the nonlinear 2-parameter model problem discussed in \cref{ssec:app:2-p}.
 The UKI-2 is applied with \dzh{$\Dt = 1/5$, $\Dt=1/3$, $\Dt =1/2$, $\Dt =2/3$, and $\Dt =3/4$}. 
Empirically we observe that all solutions considered converge to approximately the same equilibrium point;
but the convergence properties depend on $\Dt.$
The convergence of the posterior mean and covariance are reported in~\cref{fig:nonlinear-converge-gamma}. 
The best convergence rate is achieved with $\Dt$ around $1/2$ for this test. We note that the case $\Dt = 1/2$ also appears in \cite{pidstrigach2021affine} under the notion of ensemble transform Langevin dynamics.

\section{Iterative Kalman Filter}
\label{appendix:CTEKI}
\subsection{Gaussian Approximation}
Iterative Kalman filters for inverse problems are obtained by applying filtering methods, over $N = \Dt^{-1}$ steps, to the dynamical system~\cref{eq:KL_dynamics} with scaled observation error $\eta_{n+1} \sim \N(0, \Dt^{-1}\Sigma_{\eta})$. Following the discussion in Subsection~\ref{sec:KF_Intro}, the conceptual algorithm can be written as 
\begin{subequations}
\label{eq:concept2}
\begin{align}
        \mean_{n+1} &= \mean_{n} + \pCov_{n+1}^{\theta y} (\pCov_{n+1}^{y y})^{-1} (y - \py_{n+1}),  \label{eq:IKF_analysis_mean} \\
         \Cov_{n+1} &= \Cov_{n} - \pCov_{n+1}^{\theta y}(\pCov_{n+1}^{y y})^{-1} {\pCov_{n+1}^{\theta y}}{}^{T}  \label{eq:IKF_analysis_cov},
\end{align}
\end{subequations}
where 
\begin{equation}
\begin{split}
    &\py_{n+1} = \E[y_{n+1}|Y_n] = \E[\G(\theta_{n+1})|Y_n], \\
    &\pCov_{n+1}^{\theta y} = \Cov[\theta_{n+1}, y_{n+1}|Y_n] = \Cov[\theta_{n+1}, \G(\theta_{n+1})|Y_n], \\
    &\pCov_{n+1}^{y y} = \mathrm{Cov}[y_{n+1}|Y_n] = \mathrm{Cov}[\G(\theta_{n+1})|Y_n] + \Dt^{-1}\Sigma_{\eta},
\end{split}
\end{equation}
and \dzh{$Y_n:= \{y_1^\dagger, y_2^\dagger,\cdots , y_{n}^\dagger\}$}, the observation set at time $n$. Different Kalman filters~(See Subsection~\ref{sec:nonlinear_Kalman}) can be applied, which lead to the IUKF, IEnKF,  IEAKF and IETKF, all considered here.
We have the following theorem about the convergence of the conceptual algorithm:
\begin{theorem}
\label{theorem:IKF}
Consider the linear case $\G(\theta) = G \theta $. Assume that the initial mean and covariance for \eqref{eq:concept2} are the prior mean and covariance. Then
the output of the iteration \eqref{eq:concept2} at $N = \Dt ^{-1}$
equals the posterior mean and covariance.
\end{theorem}



Related calculations in continuous time may be found in Section 2.2
of \cite{garbuno2020interacting}.

\subsection{Gaussian Initialization}
\label{appendix:Gaussian_Init}
\Cref{theorem:IKF} can be extended to apply to the ensemble approximations
of the mean field dynamics \eqref{eq:concept2}, provided that the initial ensemble represents the prior mean and covariance exactly. For iterative ensemble filters, we first sample $\{\theta_j\}_{j=1}^J \sim \N(m,C)$ ($J\geq N_{\theta}+1$), and then correct them as follows. 

Define 
\begin{align}
\Theta' = [\theta_1 - m' ; \theta_2- m' ; \cdots \theta_J- m']^T \qquad m' = \frac{1}{J}\sum_{j=1}^{J} \theta_j, 
\end{align}
and correct $\Theta'$  with a $N_{\theta} \times N_{\theta}$ matrix $X$ to keep zero mean and match the covariance
\begin{align}
    \Theta := \Theta' X + 1 \otimes m^T,
\end{align}
where 
\begin{equation}
\begin{aligned}
&\textrm{SVD}: \Theta' = U_1 S_1 V_1^T, \\
&\textrm{SVD} : (J-1)C = U_2 S_2 U_2^T,\\
&X = V_1 S_1^{-1} S_2^{\frac{1}{2}} U_2.
\end{aligned}
\end{equation}

\section*{References}
\bibliographystyle{unsrt}
\bibliography{references}
\end{document}